\documentclass[reqno]{amsart}
\usepackage{amsmath}
\usepackage{amssymb}
\usepackage{hyperref}
\hypersetup{ colorlinks = true, urlcolor = blue, linkcolor = blue, citecolor = red }
\usepackage[margin=1in]{geometry}
\usepackage{xcolor}
\usepackage{comment}
\usepackage{enumitem}
\usepackage{nameref}
\DeclareMathOperator{\tr}{tr}

\DeclareMathOperator{\supp}{supp}
\DeclareMathOperator{\loc}{loc}
\DeclareMathOperator{\USC}{USC}
\DeclareMathOperator{\LSC}{LSC}
\newcommand{\R}{\mathbb{R}}
\newcommand{\ov}{\overline}
\newcommand{\norm}[1]{\left\lVert#1\right\rVert}
\newcommand{\M}{\mathcal{M}_{\lambda,\Lambda}}
\theoremstyle{plain}
\newtheorem{thm}{Theorem}[section]
\newtheorem{lem}{Lemma}[section]
\newtheorem{prop}{Proposition}[section]
\newtheorem{cor}{Corollary}[section]
\theoremstyle{definition}
\newtheorem{defn}{Definition}[section]
\theoremstyle{remark}
\newtheorem{rmk}{Remark}[section]
\numberwithin{equation}{section}
\title[]{On the fully nonlinear degenerate parabolic quenching problem}

\author[S. Kim]{Seunghyun Kim}
\address{Applied Mathematics and Computational Sciences (AMCS), Computer, Electrical and Mathematical Sciences and Engineering Division (CEMSE), King Abdullah University of Science and Technology (KAUST), Thuwal, 23955-6900, Kingdom of Saudi Arabia}
\email{seung.kim@kaust.edu.sa}

\author[J.M.~Urbano]{Jos\'{e} Miguel Urbano}
\address{Applied Mathematics and Computational Sciences (AMCS), Computer, Electrical and Mathematical Sciences and Engineering Division (CEMSE), King Abdullah University of Science and Technology (KAUST), Thuwal, 23955-6900, Kingdom of Saudi Arabia; and CMUC, Department of Mathematics, University of Coimbra, 3000-143 Coimbra, Portugal}{} 
\email{miguel.urbano@kaust.edu.sa}



\begin{document}
\everymath{\displaystyle}
\begin{abstract}
We study the fully nonlinear degenerate parabolic free
boundary problem
$$|Du|^\kappa F(D^2u)-u_t=u^{\gamma-1} \qquad \textrm{in} \quad \{u>0\},$$
where $\kappa>0$, $0<\gamma<1$, and $F$ is a convex or
concave uniformly elliptic operator with $F(0)=0$.
We establish the existence of a nonnegative viscosity
solution with prescribed continuous nonnegative boundary
data and prove its optimal upper growth along the free
boundary, together with matching non-degeneracy estimates
in parabolic cylinders adapted to the scaling.
Our approach combines a weak Bernstein estimate with
barrier arguments to obtain local compactness and boundary
attainment.
To obtain the optimal growth rate, we further use
persistence of positivity for limits with vanishing
absorption coefficients to identify the homogeneous
equation across the zero set, and then combine the
optimal spatial regularity estimate of \cite{LLYZ25}
with an iteration argument.
\end{abstract}
\maketitle

\section{Introduction}
In this paper, we study nonnegative viscosity solutions to the parabolic free boundary problem
\begin{equation}\label{eq}
\left\{\begin{aligned}
 |Du|^{\kappa}F(D^{2}u)-u_{t}
 &=u^{\gamma-1}&&\text{in }\{u>0\}\cap\Omega_{T},\\
 u&=\varphi&&\text{on }\partial_{p}\Omega_{T},
\end{aligned}\right.
\end{equation}
where $\kappa>0$, $0<\gamma<1$, and $F$ is a convex or concave uniformly elliptic operator with $F(0)=0$, as specified in \ref{a1} and \ref{a2}.
For $T>0$, let $\Omega_{T}:=\Omega\times(-T,0]$, where $\Omega\subset\R^{n}$ is a bounded domain satisfying the exterior ball condition \ref{a3}, and let $\varphi$ be a given nonnegative continuous function on the parabolic boundary $\partial_{p}\Omega_{T}$. 
The factor $|Du|^{\kappa}$ vanishes at critical points, whereas the singular absorption term $u^{\gamma-1}$ becomes unbounded as $u$ approaches zero.
The interplay between these two effects governs the growth of solutions near the free boundary $\partial\{u>0\}$.
Our aim is to establish optimal growth estimates along the free boundary and matching non-degeneracy estimates in backward parabolic cylinders adapted to the scaling.

In the stationary case with $\kappa=0$ and $F(X)=\tr X$, the equation in \eqref{eq} reduces to the positive-phase equation of the classical elliptic Alt--Phillips problem, also known as the quenching problem.
For $0<\gamma<1$, nonnegative minimizers of the functional
\begin{align*}
 \int_{\Omega}
 \left(\frac{1}{2}|Dv|^{2}+\frac{1}{\gamma}v^{\gamma}\right)\,dx
\end{align*}
with prescribed boundary values satisfy $\Delta v=v^{\gamma-1}$ in $\{v>0\}\cap\Omega$.
The regularity of minimizers and the geometry of their free boundaries were studied in \cite{D1983,D1983-2,AP1986,W1998}.
Replacing the potential $v^\gamma/\gamma$ by
$\chi_{\{v>0\}}$ yields the one-phase Bernoulli (Alt--Caffarelli) problem~\cite{AC1981}, while the potential $v$ gives the classical obstacle problem~\cite{C1977}.
Their fully nonlinear counterparts were studied in \cite{RT2011,LS2001}, respectively.

For fully nonlinear elliptic equations with singular absorption, the authors in \cite{AT2013} used a singular perturbation approach to establish optimal growth along the free boundary and non-degeneracy estimates.
A class of elliptic equations including $|Dv|^{\kappa}\Delta v=v^{\gamma-1}$ in $\{v>0\}$, with $\kappa\geq0$ and $0<\gamma<1$, was studied in \cite{T2018}.
Singular perturbation methods were also employed for fully nonlinear elliptic equations with varying exponents of gradient degeneracy and singular absorption \cite{KU26}, infinity Laplacian equations with singular absorption \cite{AS2022}, and fully nonlinear equations with logarithmic singular absorption \cite{BK26}.
In a different direction, degenerate one-phase free boundary problems were studied in \cite{DS21}, while the fully nonlinear equation $F(D^{2}v)=v^{\gamma-1}$ with $1<\gamma<2$ was investigated in \cite{WY22}.

We now turn to parabolic problems with singular absorption, which are more closely related to the present work.
Existence and quenching properties for semilinear heat equations were studied in \cite{DM05}.
Nonuniqueness of continuous weak solutions to the heat equation with singular absorption was established in \cite{W07}.
For the Cauchy problem with Laplacian diffusion, the authors in \cite{ASP26} constructed nonnegative weak solutions and established optimal growth along the free boundary and non-degeneracy estimates.
For fully nonlinear uniformly parabolic equations, Jensen--Ishii's lemma was combined with a barrier argument in \cite{ASU2024} to establish upper growth estimates along the free boundary with exponents strictly less than $2/(2-\gamma)$.
By applying Jensen--Ishii's lemma twice, the authors in \cite{KU26} obtained spatial Lipschitz estimates for a suitable power of the approximating solutions and established the optimal upper growth estimate of order $r^{2/(2-\gamma)}$ along the free boundary.

For parabolic equations with degenerate diffusion and singular absorption, the authors in \cite{DD16} studied a one-dimensional problem with the $p$-Laplacian, $p>2$.
In the positive set, the equation agrees with \eqref{eq} for $\kappa=p-2$ and $F(X)=(p-1)X$.
Using Bernstein's method after a power transformation of the regularized solution $z$, the authors obtained gradient estimates of the form
\begin{align*}
|z_x|\leq Cz^{\gamma/p},
\end{align*}
away from the initial time, with $C$ independent of the regularization parameters.

Regularity results are also available for degenerate parabolic equations with bounded right-hand side.
Under additional regularity assumptions on $F$, the authors in \cite{LLY2024} established interior H\"older estimates for the spatial gradient of viscosity solutions to
\begin{align*}
 |Du|^{\kappa}F(D^{2}u)-u_t=f
 \quad\text{in }Q_1,
\end{align*}
where $f\in C(Q_1)\cap L^{\infty}(Q_1)$.
For the homogeneous equation, local $L^{\infty}$ estimates for $u_t$ and optimal spatial $C^{1,1/(1+\kappa)}$ regularity were obtained in \cite{LLYZ25}.
We use these regularity results in our analysis, but additional estimates are needed near the zero set, where the absorption term in \eqref{eq} becomes unbounded.

We briefly describe the scaling of the equation in \eqref{eq}.
For $r>0$, define
\begin{align*}
 u_r(x,t):=\frac{u(rx,r^{\theta}t)}{r^{\alpha}},\qquad
 F_r(X):=r^{2-\alpha}F(r^{\alpha-2}X),
\end{align*}
where $\alpha,\theta>0$ are to be determined and $u_r$ is defined in $(r^{-1}\Omega)\times(-Tr^{-\theta},0]$.
Then $F_r$ satisfies \ref{a1} and \ref{a2} with the same ellipticity constants as $F$, and $u_r$ satisfies
\begin{align*}
 |Du_r|^{\kappa}F_r(D^{2}u_r)-(u_r)_t
 =u_r^{\gamma-1}
 \quad\text{in }\{u_r>0\},
\end{align*}
provided that
\begin{align*}
 \kappa(\alpha-1)+\alpha-2
 =\alpha-\theta=\alpha(\gamma-1).
\end{align*}
These relations give
\begin{align*}
 \alpha=\frac{\kappa+2}{\kappa+2-\gamma},
 \qquad
 \theta=\alpha(2-\gamma),
\end{align*}
with $1<\alpha<\theta<2$.
For $r>0$ and $z_0=(x_0,t_0)$, we define the backward parabolic cylinder
\begin{align*}
 Q_r^\theta(z_0)
 :=B_r(x_0)\times(t_0-r^\theta,t_0].
\end{align*}

We now state the main results.
The first theorem establishes the existence of a nonnegative viscosity solution and its optimal upper growth estimate along the free boundary.
\begin{thm}\label{thm1}
Let $\kappa>0$, $0<\gamma<1$, and let
$\varphi\in C(\partial_{p}\Omega_{T})$ be nonnegative. Assume
\ref{a1}--\ref{a3}, and set
\begin{align*}
 \alpha=\frac{\kappa+2}{\kappa+2-\gamma},\qquad
 \theta=\alpha(2-\gamma),\qquad
 M=\norm{\varphi}_{L^{\infty}(\partial_{p}\Omega_{T})}.
\end{align*}
There exists a viscosity solution $u\in C(\ov{\Omega_{T}})$ to
\eqref{eq} such that $0\leq u\leq M$.
If $u(z_{0})=0$ and $Q_{R}^{\theta}(z_{0})\Subset\Omega_{T}$,
then
\begin{equation}\label{growth}
 \sup_{Q_{r}^{\theta}(z_{0})}u
 \leq C\left[1+
 \left(\frac{M}{R^{\alpha}}\right)^{\frac{\kappa+2-\gamma}{2-\gamma}}
 \right]r^{\alpha} \qquad \text{for } 0<r\leq R,
\end{equation}
where $C=C(n,\lambda,\Lambda,\kappa,\gamma)$.
Moreover, if
 $B_{R}(x_{0})\times
 (t_{0}-R^{\theta},t_{0}+R^{\theta})\Subset\Omega_{T}$,
then the same bound holds with the supremum taken over
$B_{r}(x_{0})\times(t_{0}-r^{\theta},t_{0}+r^{\theta})$.
\end{thm}

We next establish a non-degeneracy estimate that follows from the viscosity subsolution inequality in the positivity set under assumption~\ref{a1} alone.
\begin{thm}\label{thm2}
Let $\kappa>0$, $0<\gamma<1$, and assume \ref{a1}. 
Let $u\in C(\Omega_{T})$ be a nonnegative viscosity subsolution to
\begin{align*}
 |Du|^{\kappa}F(D^{2}u)-u_{t}=u^{\gamma-1}
 \qquad\text{in }\{u>0\}\cap\Omega_{T}.
\end{align*}
Then there exists a constant $c_{0}=c_{0}(n,\Lambda,\kappa,\gamma)>0$ such that, for every $z_{0}\in\ov{\{u>0\}}\cap\Omega_{T}$ and
$Q_{r}^{\theta}(z_{0})\Subset\Omega_{T}$, there holds
\begin{equation}\label{nondegeneracy}
 \sup_{Q_{r}^{\theta}(z_{0})}u\geq c_{0}r^{\alpha}.
\end{equation}
\end{thm}

\begin{rmk}
Theorems~\ref{thm1} and \ref{thm2} do not require $F$ to be $1$-homogeneous.
This assumption is imposed in \cite{ASU2024,KU26}.
\end{rmk}

\begin{rmk}
The methods used to prove Theorems~\ref{thm1} and \ref{thm2} can also be adapted to the quasilinear counterpart of \eqref{eq}, where the diffusion is given by
\begin{align*}
 |Du|^{\kappa+2-p}\Delta_pu,
 \qquad \kappa>0,\quad p>1,
\end{align*}
and $\Delta_p$ denotes the $p$-Laplacian.
In particular, choosing $\kappa=p-2$ recovers the $p$-Laplacian for $p>2$.
The singular case $-1<\kappa<0$ for \eqref{eq} is the subject of ongoing work and will be treated in a forthcoming paper.
\end{rmk}

\subsection{Outline of the proof}
We construct a viscosity solution $u$ to \eqref{eq}
as a locally uniform limit of positive solutions $u_\epsilon$
to penalized problems. We define these problems by replacing
the singular absorption with the smooth truncation $g_\epsilon$
from Section~\ref{sec3} and prescribing boundary data
$\varphi+\epsilon$.

The first difficulty is to obtain estimates uniform in the
approximation parameters and strong enough to ensure local
compactness. In the related problems considered in
\cite{AT2013,ASU2024,KU26}, spatial estimates are derived
by applying Jensen--Ishii's lemma to the function obtained
by a power transformation suggested by scaling.
For \eqref{eq}, the corresponding substitution is $u=v^\alpha$.
Even for $F(X)=\tr X$, however, normalizing the coefficient
of $|Dv|^\kappa\Delta v$ leaves the factor
$\alpha^{-\kappa}v^{-\kappa(\alpha-1)}$ multiplying $v_t$,
which becomes singular as $v\to0$.
When we apply the viscosity inequalities at the two contact
points in the doubling of variables argument, the time
derivative terms are multiplied by different coefficients
and no longer cancel as in the uniformly parabolic case.
A direct adaptation of those arguments does not provide
the uniform estimates needed here.

To overcome this difficulty, we derive a weak Bernstein
estimate for smooth approximating solutions.
We first replace $F$ by its smooth approximation $F_\zeta$
and $|Du|^\kappa$ by $(\nu^2+|Du|^2)^{\kappa/2}$.
Using sub- and supersolution barriers at the parabolic
boundary, we obtain positive solutions $u_{\epsilon,\nu,\zeta}$
to the resulting problem by comparison and Perron's method.
For fixed parameters, the convexity or concavity of $F$,
together with the gradient regularity result of
\cite{LLY2024}, parabolic Schauder estimates, and a bootstrap
argument, provides the interior smoothness needed for
Bernstein's method.
For a smooth regularized solution $w$, we apply Bernstein's
method to $\eta |Dw|/w^{a/(\kappa+2)}$, where $\eta$ is a
cutoff function. Choosing $a\in(0,\gamma)$ sufficiently small
allows us to absorb the Hessian terms and control the
contribution of the singular absorption, yielding the
uniform local estimate
\begin{align*}
 \left|D\left(w^{1-\frac{a}{\kappa+2}}\right)\right|^{\kappa+2}
 =
 \left(1-\frac{a}{\kappa+2}\right)^{\kappa+2}
 \frac{|Dw|^{\kappa+2}}{w^a}
 \leq C,
\end{align*}
with $C$ independent of the regularization parameters
$\epsilon$, $\nu$, and $\zeta$.
Applying this estimate to $w=u_{\epsilon,\nu,\zeta}$,
we obtain local spatial Lipschitz estimates for
$u_{\epsilon,\nu,\zeta}^{1-a/(\kappa+2)}$,
with constants independent of $\epsilon$, $\nu$, and $\zeta$.

At fixed $\epsilon$, a barrier argument gives a time modulus
uniform in $\nu$ and $\zeta$. Together with the spatial
estimate, this yields $u_{\epsilon,\nu,\zeta}\to u_\epsilon$
locally uniformly as $(\nu,\zeta)\to(0,0)$.
To obtain time estimates uniform in $\epsilon$, we use the
power $u_\epsilon^{2-a}$. In its viscosity inequalities,
the absorption is controlled by $u_\epsilon^{\gamma-a}$,
while the additional gradient term is controlled by the
Bernstein estimate. Comparison with functions quadratic
in space and affine in time then gives a uniform time
modulus. Together with the spatial estimate, this yields
local compactness as $\epsilon\to0$.

The boundary barriers used to construct $u_{\epsilon,\nu,\zeta}$
also give estimates near the lateral and initial boundaries
that are uniform in $\epsilon$, $\nu$, and $\zeta$.
The lower barriers use power profiles $W$ with exponent
$\alpha$, which balances gradient degeneracy and singular
absorption. These profiles satisfy
\begin{align*}
 (\nu^2+|DW|^2)^{\kappa/2}F_\zeta(D^2W)-W_t
 \geq W^{\gamma-1}
 \geq g_\epsilon(W)
 \quad\text{in }\{W>0\}\cap\Omega_T.
\end{align*}
The resulting boundary estimates ensure that
$u_\epsilon=\varphi+\epsilon$ on $\partial_p\Omega_T$
after passing to the limit as $(\nu,\zeta)\to(0,0)$
with $\epsilon$ fixed, and that
$u=\varphi$ on $\partial_p\Omega_T$
after letting $\epsilon\to0$.
Together with local compactness and monotonicity in
$\epsilon$, these estimates yield a viscosity solution
$u\in C(\overline{\Omega_T})$ to \eqref{eq}.

The Bernstein estimate gives spatial growth at zeros with
exponent $(\kappa+2)/(\kappa+2-a)$, which is strictly smaller
than $\alpha$.
To recover the optimal exponent, we consider normalized
positive solutions $v_j$ with absorption terms
$\delta_jg_{\epsilon_j}(v_j)$, where $\delta_j\to0$, and
pass to a locally uniform limit $V$.
The absorption tends to zero locally in $\{V>0\}$, but
need not do so uniformly near $\{V=0\}$.
Stability therefore identifies the homogeneous equation
directly only in the positivity set.
The Hopf-type argument used in \cite{AT2013,KU26}
for elliptic problems does not directly apply, since
spatial slices of the limiting parabolic solution need
not satisfy a homogeneous elliptic equation.

The key step is to identify the homogeneous equation
across the zero set.
For sufficiently large $j$, we compare $v_j$ with the
positive subsolution
\begin{align*}
A(t)\phi\left(\frac{x-x_0}{r}\right)
\end{align*}
constructed in Lemma~\ref{comparison bump}.
Passing to the limit shows that positivity of $V$ at a
fixed spatial point persists forward in time within the
cylinder (Lemma~\ref{positivity}).
Consequently, if $V(x_0,t_0)=0$, then $V(x_0,t)=0$ at
preceding times.
Any test function touching $V$ from below at $(x_0,t_0)$
therefore has nonnegative time derivative there.
Its spatial gradient vanishes by the superlinear spatial
growth supplied by the weak Bernstein estimate.
Since $\kappa>0$, these two facts give the homogeneous
supersolution inequality at the zero point.
The subsolution inequality already follows throughout
the cylinder from the sign of the absorption and the
stability of viscosity solutions.
Thus, for a limiting operator $F_\infty$ satisfying
\ref{a1} and \ref{a2}, we obtain
\begin{align*}
 |DV|^\kappa F_\infty(D^2V)-V_t=0
 \quad\text{throughout the cylinder}.
\end{align*}

Persistence of positivity has a second role in the
improvement argument.
If $V(0,0)=0$, it gives $V(0,t)=0$ on the preceding
time interval.
Applying the optimal spatial regularity estimate of
\cite{LLYZ25} on each time slice then gives
\begin{align*}
 V(x,t)\leq C|x|^{(\kappa+2)/(\kappa+1)}
\end{align*}
near the origin.
Here we use the convexity or concavity of $F_\infty$,
inherited from assumption~\ref{a2}.
Since $(\kappa+2)/(\kappa+1)>\alpha$, compactness yields
a one-step improvement at a fixed scale.
Iteration gives the optimal upper growth estimate for $u_\epsilon$ in backward cylinders. 
We extend the resulting spatial bound on the central time slice to later times by comparison with a supersolution that is radial in space.
Passing to the limit as $\epsilon\to0$ proves Theorem~\ref{thm1}.

Finally, we establish the non-degeneracy estimate directly
from the viscosity subsolution inequality under
assumption~\ref{a1} alone.
For $z_0=(x_0,t_0)\in\{u>0\}\cap\Omega_T$, set
$M_r=\sup_{Q_r^\theta(z_0)}u$.
We use a comparison function quadratic in space and affine
in time, adding the temporal penalization $\eta/(t_0-t)$
to obtain an interior upper contact point in the positivity
set before $t_0$.
The viscosity subsolution inequality at this point gives
\begin{align*}
 1\leq C\left(\frac{M_r}{r^\alpha}\right)^{\kappa+2-\gamma}
       +\left(\frac{M_r}{r^\alpha}\right)^{2-\gamma}.
\end{align*}
Consequently, $M_r\geq c_0r^\alpha$.
By continuity, \eqref{nondegeneracy} also holds for
$z_0\in\overline{\{u>0\}}\cap\Omega_T$.
The argument applies to any nonnegative viscosity subsolution
and does not require the approximation procedure.

\subsection{Organization of the paper}
Section~\ref{sec2} introduces the notation, assumptions, and
preliminary results.
Section~\ref{sec3} constructs smooth regularized solutions
and establishes boundary estimates uniform in the
regularization parameters.
In Section~\ref{sec4}, a weak Bernstein estimate and
temporal equicontinuity estimates yield local compactness
and the existence of a viscosity solution to \eqref{eq}.
In Section~\ref{sec5}, we show that as the absorption coefficient tends to zero, each limiting function satisfies the homogeneous equation across its zero set. 
Compactness and iteration give the optimal upper growth estimate in backward cylinders.
We then extend this estimate forward in time by comparison.
Finally, Section~\ref{sec6} establishes the non-degeneracy estimate using the viscosity subsolution inequality.

\section{Preliminaries}
\label{sec2}

\subsection{Notation, definitions and main assumptions}
Let $\Omega\subset\R^{n}$ be a bounded domain and $T>0$.
We denote the space-time cylinder and its parabolic boundary by
\begin{align*}
 \Omega_T&:=\Omega\times(-T,0],\\
 \partial_p\Omega_T
 &:=(\overline{\Omega}\times\{-T\})
   \cup(\partial\Omega\times[-T,0]).
\end{align*}
For $z_0=(x_0,t_0)$, $r>0$, and $\theta>0$,
we denote the open ball and backward parabolic cylinder by
\begin{align*}
 B_r(x_0)&:=\{x\in\R^n:|x-x_0|<r\},\\
 Q_r^\theta(z_0)&:=B_r(x_0)\times(t_0-r^\theta,t_0].
\end{align*}
When the center is the origin, we omit it.
We also write $Q_{r}=Q_{r}^{2}$.
The parabolic distance is defined by
\begin{align*}
 d_{p}((x,t),(y,s)):=|x-y|+|t-s|^{1/2}.
\end{align*}
For $\sigma\in(0,1)$, $C^{\sigma,\sigma/2}$ and
$C^{2+\sigma,1+\sigma/2}$ denote the standard parabolic
H\"older spaces.

 Let $\mathcal{S}^{n}$ denote the space of real $n\times n$ symmetric matrices. 
 For a matrix $X=(X_{ij})$, we write $|X|:=(\sum_{i,j}X_{ij}^{2})^{1/2}$ for its Frobenius norm.
 For parameters $0<\lambda\leq\Lambda$, the Pucci extremal operators $\mathcal{M}^{\pm}_{\lambda, \Lambda}: \mathcal{S}^{n}\rightarrow \R$ are defined as 
\begin{align*}
\mathcal{M}^{+}_{\lambda, \Lambda}(M):=  \Lambda\sum_{e_i>0}e_{i}+\lambda\sum_{e_{i}<0}e_{i} \quad \text{and} \quad \mathcal{M}^{-}_{\lambda, \Lambda}(M):=  \lambda\sum_{e_{i}>0}e_{i}+\Lambda\sum_{e_{i}<0}e_{i}, 
\end{align*}
where $e_{i}=e_{i}(M)$ are the eigenvalues of $M\in\mathcal{S}^{n}$. 

We recall the definition of viscosity solutions; see \cite[Section~8]{CIL1992}.
\begin{defn}
\label{viscosity solution}
Let $g\in C(\R\times\Omega_{T})$. An upper semicontinuous function
$u:\Omega_{T}\to\R$ is a viscosity subsolution to
\begin{align}
\label{def-eq}
 |Du|^{\kappa}F(D^{2}u)-u_{t}=g(u,x,t)
 \quad\text{in }\Omega_{T}
\end{align}
if, whenever a test function $\phi\in C^{2,1}$ touches $u$ from above at $z_{0}=(x_{0},t_{0})\in\Omega_{T}$, one has
\begin{align*}
 |D\phi(z_{0})|^{\kappa}F(D^{2}\phi(z_{0}))
 -\phi_{t}(z_{0})\geq g(u(z_{0}),z_{0}).
\end{align*}
A lower semicontinuous function $u:\Omega_{T}\to\R$ is a viscosity
supersolution if the inequality is reversed for test functions
touching from below.
A continuous function is a viscosity solution to \eqref{def-eq}
if it is both a viscosity subsolution and a viscosity supersolution.
At the terminal time $t_{0}=0$, the touching conditions are
understood relative to $\Omega_{T}$.
\end{defn}

We now state the main assumptions in this paper.
\begin{enumerate}[label=(A\arabic*),ref=(A\arabic*)]
 \item\label{a1} 
 $F(0)=0$ and $F$ is $(\lambda, \Lambda)$-uniformly elliptic, \textit{i.e.}, 
$$\mathcal{M}^{-}_{\lambda, \Lambda}(X-Y)\leq F(X)-F(Y)\leq\mathcal{M}^{+}_{\lambda, \Lambda}(X-Y),$$ 
for every $X,Y\in \mathcal{S}^{n}$.
 \item\label{a2} $F$ is convex or concave.
 \item\label{a3} $\Omega$ satisfies the exterior ball condition:
for every $x_{0}\in\partial\Omega$, there exist $r_{0}>0$ and
$y_{0}\in\R^{n}\setminus\overline{\Omega}$ such that
$\overline{B_{r_{0}}(y_{0})}\cap\overline{\Omega}=\{x_{0}\}.$
\end{enumerate}

\subsection{Preliminary results}

We begin with an elementary observation that allows us
to use spatial Lipschitz bounds in viscosity inequalities.

\begin{lem}
\label{lipschitz test gradient}
Let $Q=U\times I$, where $U\subset\R^{n}$ is open and
$I\subset\R$ is an interval.
Suppose that $w\in C(Q)$ satisfies
\begin{align*}
|w(x,t)-w(y,t)|\leq L|x-y|
\qquad\text{for }(x,t),(y,t)\in Q,
\end{align*}
where $L\geq0$.
If $\xi\in C^{2,1}$ touches $w$ from above or below at
$z_{0}\in Q$, then $|D\xi(z_{0})|\leq L$.
\end{lem}

\begin{proof}
Assume that $\xi$ touches $w$ from above at
$z_{0}=(x_{0},t_{0})$.
For every $|e|=1$ and all sufficiently small $h>0$, we have
\begin{align*}
-Lh
\leq w(x_{0}+he,t_{0})-w(z_{0})
\leq \xi(x_{0}+he,t_{0})-\xi(z_{0}).
\end{align*}
Dividing by $h$ and letting $h\to0^{+}$, we obtain
\begin{align*}
D\xi(z_{0})\cdot e\geq-L
\qquad\text{for every }|e|=1.
\end{align*}
If $D\xi(z_{0})\neq0$, taking
$e=-D\xi(z_{0})/|D\xi(z_{0})|$ yields
$-|D\xi(z_{0})|\geq-L$, which implies $|D\xi(z_{0})|\leq L$.
The same bound is immediate if $D\xi(z_{0})=0$.
The case of touching from below follows by applying the same argument to $-w$ and $-\xi$.
\end{proof}

The next two results provide the regularity needed to obtain
smooth approximating solutions in Section~\ref{sec3}.
We first recall an interior gradient regularity estimate.
\begin{lem}[{\cite[Theorem 1.1]{LLY2024}}]
\label{gradient-regularity}
Let $F\in C^{\infty}(\mathcal{S}^{n})$ satisfy \ref{a1} and
\ref{a2}, and suppose that $DF$ is globally Lipschitz.
Let $\kappa>0$ and $f\in C(Q_{1})\cap L^{\infty}(Q_{1})$.
Let $u$ be a bounded viscosity solution to
\begin{align*}
 |Du|^{\kappa}F(D^{2}u)-u_{t}=f
 \quad\text{in }Q_{1}.
\end{align*}
Then, there exists a constant $\sigma\in(0,1)$ depending only on $n, \lambda, \Lambda, \kappa, \norm{u}_{L^{\infty}(Q_{1})}$, and $\norm{F}_{C^{1,1}(\mathcal{S}^{n})}$ such that $u\in C^{1+\sigma, (1+\sigma)/2}(\ov{Q_{1/2}})$ with an estimate
\begin{align*}
\norm{u}_{C^{1+\sigma,(1+\sigma)/2}(Q_{1/2})}\leq C,
\end{align*}
where $C>0$ is a constant depending only on $n, \lambda, \Lambda, \kappa, \norm{u}_{L^{\infty}(Q_{1})}$, and $\norm{f}_{L^{\infty}(Q_{1})}$.
\end{lem}

We also use the following consequence of the parabolic Schauder
estimates.
\begin{lem}[{\cite[Theorem 1.1]{W1992}}]
\label{schauder}
Let $G(x,t,X)$ be uniformly elliptic and either convex or concave in $X$. 
Assume that $f\in C^{\sigma,\sigma/2}_{\loc}(\Omega_{T})$ for some $\sigma\in(0,1)$.
In addition, assume that for every cylinder $Q\Subset\Omega_{T}$, there exists $C_{Q}>0$ such that
\begin{align}
\label{coeff holder}
 \sup_{X\in\mathcal{S}^{n}}
 \frac{|G(x,t,X)-G(y,s,X)|}{1+|X|}
 \leq C_{Q}d_{p}((x,t),(y,s))^{\sigma}
\end{align}
for all $(x,t),(y,s)\in Q$.
Let $u$ be a bounded viscosity solution to
\begin{align*}
 u_{t}=G(x,t,D^{2}u)+f(x,t)
 \quad \text{in } \Omega_{T}.
\end{align*}
Then $u\in C^{2+\eta,1+\eta/2}_{\loc}(\Omega_{T})$ for some $\eta\in(0,\sigma]$.
\end{lem}

The following optimal spatial gradient estimate for the
homogeneous equation will be used in the improvement
argument of Section~\ref{sec5}.
\begin{lem}[{\cite[Theorem 1.3]{LLYZ25}}]
\label{time derivative}
Let $F$ satisfy \ref{a1} and \ref{a2}, and let $\kappa>0$.
Let $v$ be a bounded viscosity solution to
\begin{align*}
 |Dv|^{\kappa}F(D^{2}v)-v_{t}=0
 \quad\text{in }Q_{1}.
\end{align*}
Then $v(\cdot,t)\in C^{1,1/(1+\kappa)}(B_{1/2})$
for every $t\in[-1/4,0]$, and
\begin{align*}
|Dv(x,t)-Dv(y,t)|
\leq C|x-y|^{1/(1+\kappa)}
\end{align*}
for all $x,y\in B_{1/2}$ and $t\in[-1/4,0]$, where
$C=C(n,\lambda,\Lambda,\kappa,\|v\|_{L^{\infty}(Q_{1})})>0$.
\end{lem}

Finally, we recall the eigenvalues of the Hessian of a radial function, which will be used in the construction of barriers.

\begin{lem}\label{Hessian}
Let $x_{0}\in\R^{n}$ and $g\in C^{2}(0,\infty)$.
For $f(x)=g(r)$, where $r=|x-x_{0}|>0$,
the eigenvalues of $D^{2}f(x)$ are $g''(r)$
in the radial direction and $g'(r)/r$ in
the $n-1$ tangential directions.
\end{lem}
\begin{proof}
Writing $e=(x-x_{0})/r$, we have
\begin{align*}
D^{2}f(x)
=g''(r)e\otimes e
+\frac{g'(r)}{r}(I-e\otimes e),
\end{align*}
which yields the desired eigenvalues.
\end{proof}

\section{Approximation}
\label{sec3}
We construct our solution to \eqref{eq} as the limit of
solutions to singularly penalized approximating problems.
To obtain the estimates needed for this limit, we further
regularize the diffusion and construct smooth solutions
to the resulting Dirichlet problems.
The boundary estimates obtained here will also be used
when the regularization parameters tend to zero.

\subsection{Approximation procedure}
We first truncate the singular absorption term.
Let $\beta:\R\to\R$ be a smooth nondecreasing function such that
\begin{align*}
&0\leq \beta \leq 1, \nonumber \\
&\beta(s)=0 \quad \textrm{for} \quad s\leq 1, \nonumber \\ 
&\beta(s)=1 \quad \textrm{for} \quad s\geq 2.
\end{align*}
For $\epsilon>0$, define
\begin{align*}
 \beta_{\epsilon}(s):=\beta(s/\epsilon),\qquad
 g_{\epsilon}(s):=
 \begin{cases}
  \beta_{\epsilon}(s)s^{\gamma-1},&s>0,\\
  0,&s\leq0.
 \end{cases}
\end{align*}
As $\epsilon \to 0$, $\beta_{\epsilon}(s)$ converges to $\chi_{\{s>0\}}$. Here, $\chi_{E}$ denotes the characteristic function of a set $E$. For every fixed $\epsilon>0$, the function $g_{\epsilon}$ is smooth,
bounded, globally Lipschitz, and nonnegative. 
We first consider the penalized problem
\begin{equation}
\label{pe1}
\left\{\begin{aligned}
 |Du_{\epsilon}|^{\kappa}F(D^{2}u_{\epsilon})
 -(u_{\epsilon})_{t}
 &=g_{\epsilon}(u_{\epsilon})
 &&\text{in }\Omega_{T},\\
 u_{\epsilon}
 &=\varphi+\epsilon
 &&\text{on }\partial_{p}\Omega_{T}.
\end{aligned}\right.
\end{equation}
To obtain the smooth approximating solutions needed for
Bernstein's method, we further regularize the gradient
factor and the operator $F$. We extend $F$ from $\mathcal{S}^{n}$ to
$\mathbb{R}^{n\times n}$ by setting
\begin{align*}
 F(M)=F\left(\frac{M+M^{T}}{2}\right),
 \qquad M\in\mathbb{R}^{n\times n},
\end{align*}
using the same notation for the extension.
Let $\psi\in C_c^{\infty}(\mathbb{R}^{n\times n})$
be a standard mollifier satisfying
\begin{align*}
 \int_{\mathbb{R}^{n\times n}}\psi(M)\,dM=1,
 \qquad
 \operatorname{supp}\psi
 \subset
 \left\{
 M\in\mathbb{R}^{n\times n}:
 |M| \leq1
 \right\}.
\end{align*}
For $\zeta>0$, set
$\psi_{\zeta}(M):=\zeta^{-n^{2}}\psi(M/\zeta)$ and define
\begin{align*}
 F_{\zeta}(X)
 :=(F*\psi_{\zeta})(X)-(F*\psi_{\zeta})(0)=\int_{\mathbb{R}^{n\times n}}
   \bigl(F(X-M)-F(-M)\bigr)\psi_{\zeta}(M)\,dM,
 \qquad X\in\mathcal{S}^{n}.
\end{align*}
Then $F_{\zeta}\in C^{\infty}(\mathcal{S}^{n})$ satisfies \ref{a1}. Moreover, $F_{\zeta}$ inherits the convexity or concavity of $F$.
Finally, since $F$ is Lipschitz continuous and $F(0)=0$, $F_{\zeta}\to F$ uniformly on $\mathcal{S}^{n}$ as $\zeta\to0$.
For $\nu\geq0$, set
\begin{align*}
 a_{\nu}(q):=(\nu^{2}+|q|^{2})^{\kappa/2},
 \qquad q\in\mathbb{R}^{n}.
\end{align*}
For $\epsilon,\nu,\zeta\in(0,1)$, we consider the
regularized problem
\begin{equation}
\label{pe}
\left\{\begin{aligned}
 a_{\nu}(Du_{\epsilon,\nu,\zeta})
 F_{\zeta}(D^{2}u_{\epsilon,\nu,\zeta})
 -(u_{\epsilon,\nu,\zeta})_{t}
 &=g_{\epsilon}(u_{\epsilon,\nu,\zeta})
 &&\text{in }\Omega_{T},\\
 u_{\epsilon,\nu,\zeta}
 &=\varphi+\epsilon
 &&\text{on }\partial_{p}\Omega_{T}.
\end{aligned}\right.
\end{equation}
The existence and interior smoothness of solutions to
\eqref{pe} will be established in Theorem~\ref{existence}.
For each fixed $\epsilon>0$, the solutions
$u_{\epsilon,\nu,\zeta}$ converge locally uniformly to the
unique viscosity solution $u_\epsilon$ to \eqref{pe1} as
$(\nu,\zeta)\to(0,0)$; see
Proposition~\ref{first compactness}.
The same construction applies with $g_\epsilon$ replaced
by $\delta g_\epsilon$ for $0\leq\delta\leq1$.
We present the construction for $\delta=1$ and use this
extension in Section~\ref{sec5}.

\subsection{Comparison principle}
We establish a comparison principle for equations with
a globally Lipschitz nonlinearity $f(u)$.
It applies to both \eqref{pe} and \eqref{pe1} and will be
used in Perron's method and in the subsequent limiting
arguments.
\begin{lem}
\label{comparison}
Let $H:\R^{n}\times\mathcal{S}^{n}\to\R$ be continuous and degenerate elliptic, i.e., $H(p,X)\leq H(p,Y)$ for every $p\in \R^{n}$ and $X,Y\in\mathcal{S}^{n}$ with $X\leq Y$.
Let $f:\R\to\R$ be globally Lipschitz, and let $v\in\USC(\ov{\Omega_{T}})$ and $w\in\LSC(\ov{\Omega_{T}})$ be viscosity sub- and supersolutions, respectively, to
\begin{align*}
 H(Du,D^{2}u)-u_{t}=f(u)
 \quad\text{in }\Omega\times(-T,0).
\end{align*}
If $v\leq w$ on $\partial_{p}\Omega_{T}$, then $v\leq w$ in
$\Omega\times(-T,0)$.
\end{lem}
\begin{proof}
Let $L=[f]_{\operatorname{Lip}(\R)}$ and choose $\tau>0$ such
that $2L\tau<1$. It suffices to prove comparison in
$\Omega\times[s,b)$ for $-T\leq s<b\leq0$ and $b-s\leq\tau$,
assuming $v(\cdot,s)\leq w(\cdot,s)$ in $\overline\Omega$.
Starting with $s=-T$ and repeating the argument on overlapping
intervals of length at most $\tau$ then proves the assertion.
Since $v$ is bounded above and $w$ is bounded below on the compact
set $\ov{\Omega_{T}}$, there is a constant $K>0$ such that
$v(x,t)-w(y,t)\leq K$ for all $x,y\in\overline\Omega$ and
$t\in[-T,0]$.
Suppose, to the contrary, that
\begin{align*}
m:=\sup_{\overline{\Omega}\times[s,b)}(v-w)>0.
\end{align*}
Choose $(\widehat x,\widehat t)\in\overline\Omega\times[s,b)$ such
that $(v-w)(\widehat x,\widehat t)>3m/4$, and fix
$0<\eta<m(b-\widehat t)/8$. 
Then we have
\begin{align*}
\sup_{\overline{\Omega}\times[s,b)}
 \left(v-w-\frac{m}{2\tau}(t-s)-\frac{\eta}{b-t}\right)
 \geq (v-w)(\widehat x, \widehat t)-\frac{m}{2\tau}(\widehat t-s)-\frac{\eta}{b-\widehat t}   >\frac{m}{8}.
\end{align*}
For $h>0$, define
\begin{align*}
 \Phi_{h}(x,y,t):=v(x,t)-w(y,t)-\frac{|x-y|^{2}}{2h}
 -\frac{m}{2\tau}(t-s)-\frac{\eta}{b-t}
\end{align*}
on $\overline\Omega\times\overline\Omega\times[s,b)$.
The function $\Phi_h$ is upper semicontinuous and $\Phi_{h}\to-\infty$
uniformly in $x,y$ as $t\to b^{-}$. Hence it attains its maximum
at some $(x_h,y_h,t_h)$, where
\begin{align}
\label{bound Phi}
 \Phi_h(x_h,y_h,t_h)>\frac{m}{8},\qquad
 \frac{|x_h-y_h|^2}{2h}+\frac{\eta}{b-t_h}\leq K.
\end{align}
In particular, $|x_h-y_h|\leq \sqrt{2Kh}\to0$ as $h\to0^{+}$ and $t_h\leq b-\eta/K<b$.
By compactness and semicontinuity,
\begin{align}
\label{diff bound}
 \limsup_{h\to0^{+}}
 \bigl(v(x_{h},t_{h})-w(y_{h},t_{h})\bigr)\leq m.
\end{align}
Moreover, $x_h,y_h\in\Omega$ and $s<t_h<b$ for all sufficiently
small $h$. 
Otherwise, there exists 
$h_j\to0^{+}$ and
$(\ov{x},\ov{t})\in\overline\Omega\times[s,b)$
such that $(x_{h_j},y_{h_j},t_{h_j})\to(\ov{x},\ov{x},\ov{t})$,
where $\ov{x}\in\partial\Omega$ or $\ov{t}=s$.
By \eqref{bound Phi}, semicontinuity, and the boundary ordering,
\begin{align*}
m/8\leq\limsup_{j\to \infty}\Phi_{h_{j}}(x_{h_{j}},y_{h_{j}},t_{h_{j}})\leq\limsup_{j\to \infty}\bigl(v(x_{h_{j}},t_{h_{j}})-w(y_{h_{j}},t_{h_{j}})\bigr)
\leq v(\ov{x},\ov{t})-w(\ov{x},\ov{t})\leq0,    
\end{align*}
which yields a contradiction.

By the parabolic theorem of sums
\cite[Theorem 8.3]{CIL1992}, there exist $a_h\in\R$ and
$X_h,Y_h\in\mathcal S^n$ with $X_h\leq Y_h$ such that
\begin{align*}
 \left(a_h+\frac{m}{2\tau}+\frac{\eta}{(b-t_h)^2},
       \frac{x_h-y_h}{h},X_h\right)
 \in\overline{\mathcal P}^{2,+}v(x_h,t_h),\quad
 \left(a_h,\frac{x_h-y_h}{h},Y_h\right)
 \in\overline{\mathcal P}^{2,-}w(y_h,t_h).
\end{align*}
The required local upper bounds on temporal superjet components
for $v$ and $-w$, when the function values and spatial jet
components are bounded, follow from their viscosity inequalities
and the continuity of $H$ and $f$.
Subtracting the viscosity inequalities, we obtain
\begin{align*}
 \frac{m}{2\tau}\leq\frac{m}{2\tau}+\frac{\eta}{(b-t_{h})^{2}}
 &\leq H\left(\frac{x_{h}-y_{h}}{h},X_{h}\right)-H\left(\frac{x_{h}-y_{h}}{h},Y_{h}\right)
       +f(w(y_{h},t_{h}))-f(v(x_{h},t_{h}))\\
 &\leq L\bigl(v(x_{h},t_{h})-w(y_{h},t_{h})\bigr).
\end{align*}
Here we used $X_h\leq Y_h$, the Lipschitz continuity of $f$, and
$v(x_h,t_h)>w(y_h,t_h)$. Letting $h\to0^{+}$ and using \eqref{diff bound}, we obtain $m/(2\tau)\leq Lm$, contradicting $2L\tau<1$. This completes the proof.
\end{proof}

\subsection{Boundary barriers and uniform boundary attainment}

We construct sub- and supersolution barriers at the
parabolic boundary with local moduli of continuity
independent of $\epsilon$, $\nu$, and $\zeta$.
These barriers will be used in Perron's method and will
ensure boundary attainment as $(\nu,\zeta)\to(0,0)$ and
then $\epsilon\to0$.
We begin with an estimate for a radial profile, which
will also be used in Section~\ref{sec5} to prove persistence
of positivity for limits with vanishing absorption
coefficients.

\begin{lem}
\label{barrier homo}
Let $\kappa>0$, $0<\gamma<1$, and set
\begin{align*}
\alpha:=\frac{\kappa+2}{\kappa+2-\gamma}, \qquad \phi(y):=(1-|y|^{2})^{\alpha} \quad \text{for }y\in\overline{B_{1}}.    
\end{align*}
There exist constants $K_{0}\geq1$ and $k_{0}>0$, depending only on
$n,\lambda,\Lambda,\kappa$, and $\gamma$, such that, for every $\nu\in[0,1]$,
\begin{align}
\label{bump-ineq}
a_{\nu}(D\phi)\M^{-}(D^{2}\phi)+K_{0}\phi
\geq k_{0}\phi^{\gamma-1}\qquad\text{in }B_{1}.
\end{align}
\end{lem}

\begin{proof}
For $y\in B_{1}$, direct differentiation gives
\begin{align*}
D\phi(y)&=-2\alpha y(1-|y|^{2})^{\alpha-1},\\
D^{2}\phi(y)&=4\alpha(\alpha-1)(1-|y|^{2})^{\alpha-2}y\otimes y
-2\alpha(1-|y|^{2})^{\alpha-1}I.
\end{align*}
Consequently,
\begin{align*}
\M^{-}(D^{2}\phi(y))
&\geq 4\lambda\alpha(\alpha-1)|y|^{2}(1-|y|^{2})^{\alpha-2}
-2n\Lambda\alpha(1-|y|^{2})^{\alpha-1}.
\end{align*}
Choose $c\in(0,1/2)$ sufficiently small that
\begin{align*}
n\Lambda(2c-c^{2})\leq\lambda(\alpha-1)(1-c)^{2}.
\end{align*}
In the region $\{1-c\leq|y|<1\}$, we have
\begin{align*}
\M^{-}(D^{2}\phi(y))
&\geq2\lambda\alpha(\alpha-1)|y|^{2}(1-|y|^{2})^{\alpha-2},\\
a_{\nu}(D\phi(y))\M^{-}(D^{2}\phi(y))&\geq|D\phi(y)|^{\kappa}\M^{-}(D^{2}\phi(y)) \\
&\geq2^{\kappa+1}\lambda\alpha^{\kappa+1}(\alpha-1)
(1-c)^{\kappa+2}(1-|y|^{2})^{\alpha(\gamma-1)},
\end{align*}
where we have used $\kappa(\alpha-1)+\alpha-2=\alpha(\gamma-1)$.
Let
\begin{align*}
k_{0}:=2^{\kappa+1}\lambda\alpha^{\kappa+1}(\alpha-1)(1-c)^{\kappa+2}.
\end{align*}
On $\overline{B_{1-c}}$, the function $\phi$ is positive and smooth.
We can therefore choose $K_{0}\geq1$ such that
\begin{align*}
K_{0}\geq
\max_{\substack{y\in\overline{B_{1-c}}\\0\leq\nu\leq1}}
\frac{k_{0}\phi(y)^{\gamma-1}-a_{\nu}(D\phi(y))\M^{-}(D^{2}\phi(y))}{\phi(y)}.
\end{align*}
This proves \eqref{bump-ineq}.
\end{proof}

We now construct supersolution barriers at the lateral
and initial boundaries.

\begin{lem}
\label{bdrysuper}
Let $\kappa>0$, $0<\gamma<1$, and assume \ref{a1} and \ref{a3}.
Let $\varphi\in C(\partial_{p}\Omega_{T})$ be nonnegative, and set
\begin{align*}
M_{0}:=\norm{\varphi}_{L^{\infty}(\partial_{p}\Omega_{T})}+1.
\end{align*}
For every $z_{0}\in\partial_{p}\Omega_{T}$, $\sigma>0$, and
$\epsilon,\nu,\zeta\in(0,1)$, there exists a viscosity supersolution
$\overline V\in C(\overline{\Omega_{T}})$ to
\eqref{pe} such that
\begin{align}
\label{super es}
\epsilon&\leq\overline V\leq M_{0}
\quad\text{in }\overline{\Omega_{T}}, \qquad
\varphi+\epsilon\leq\overline V
\quad\text{on }\partial_{p}\Omega_{T},\qquad \varphi(z_{0})+\epsilon\leq\overline V(z_{0})
\leq\varphi(z_{0})+\epsilon+\sigma.
\end{align}
Moreover, $\overline V$ has a modulus of continuity independent of $\epsilon,\nu,\zeta$ in a fixed relative neighborhood of $z_{0}$.
\end{lem}

\begin{proof}
Fix $z_{0}=(x_{0},t_{0})\in\partial_{p}\Omega_{T}$ and $\sigma>0$,
and write $m:=\varphi(z_{0})+\epsilon$.
We divide the proof into two cases according to the location of $z_{0}$.

\smallskip
\noindent\textit{Case I. The lateral boundary.}
Suppose that $x_{0}\in\partial\Omega$ and $-T<t_{0}\leq0$.
Let $B_{2\rho}(y_{0})$ be an exterior ball tangent at $x_{0}$, and set
\begin{align*}
y:=\frac{x_{0}+y_{0}}{2},\qquad
\eta(x):=|x-y|-\rho,\qquad
\ell:=\max\left\{1,\frac{2\Lambda(n-1)}{\lambda\rho}\right\}.
\end{align*}
By continuity of $\varphi$ at $z_{0}$, choose $r>0$ and
$0<h<t_{0}+T$ sufficiently small that
\begin{align*}
|\varphi(z)-\varphi(z_{0})|<\sigma/2
\quad\text{for }z\in\partial_{p}\Omega_{T}
\cap\bigl(\overline{B_{r}(x_{0})}\times[t_{0}-h,t_{+}]\big),
\end{align*}
where $t_{+}:=\min\{0,t_{0}+h\}$.
We consider the barrier
\begin{align*}
\overline W(x,t):=m+\sigma+A(1-e^{-\ell\eta(x)})
+\frac{M_{0}}{h^{2}}(t-t_{0})^{2}
\quad\text{for }(x,t)\in\overline{\Omega_{T}},
\end{align*}
where $A>0$ is to be determined.
The exterior ball condition gives
\begin{align*}
\eta(x)\geq\frac{|x-x_{0}|^{2}}{2(|x-y|+\rho)}
\geq\frac{|x-x_{0}|^{2}}{2(|x-x_{0}|+2\rho)}
\quad\text{for }x\in\overline\Omega,
\end{align*}
and hence
\begin{align}
\label{eta bound}
0\leq\eta\leq\operatorname{diam}\Omega\quad\text{in }\overline\Omega,\qquad
\eta\geq\frac{r^{2}}{2(r+2\rho)}=:\mu>0
\quad\text{on }\overline\Omega\setminus B_{r}(x_{0}).
\end{align}
In $\Omega_{T}$, a direct calculation gives $|D\overline W|=A\ell e^{-\ell\eta}$ and
\begin{align*}
F_{\zeta}(D^{2}\overline W)\leq\M^{+}(D^{2}\ov{W})\leq A\ell e^{-\ell\eta}
\left(-\lambda\ell+\frac{\Lambda(n-1)}{|x-y|}\right)
\leq-\frac{\lambda}{2}A\ell^{2}e^{-\ell\eta}<0,
\end{align*}
where we used \ref{a1}, Lemma~\ref{Hessian}, and the choice of $\ell$.
Choose $A$ large enough that
\begin{align*}
A(1-e^{-\ell\mu})\geq M_{0},\qquad
\frac{\lambda}{2}\ell^{\kappa+2}
e^{-(\kappa+1)\ell\operatorname{diam}\Omega}A^{\kappa+1}\geq\frac{2M_{0}T}{h^{2}}.
\end{align*}
Since $a_{\nu}(D\overline W)\geq|D\overline W|^{\kappa}$,
in $\Omega_{T}$ we obtain
\begin{align*}
a_{\nu}(D\overline W)F_{\zeta}(D^{2}\overline W)-(\overline W)_{t}
\leq-\frac{\lambda}{2}\ell^{\kappa+2}
e^{-(\kappa+1)\ell\operatorname{diam}\Omega}A^{\kappa+1}+\frac{2M_{0}T}{h^{2}}\leq0\leq g_{\epsilon}(\ov{W}).
\end{align*}
Thus $\overline W$ is a viscosity supersolution in $\Omega_{T}$.
Furthermore,
\begin{align*}
\ov{W}&\geq m+\sigma/2\geq \varphi+\epsilon
\quad\text{on } \partial_{p}\Omega_{T}
\cap\bigl(\overline{B_{r}(x_{0})}\times[t_{0}-h,t_{+}]\big),\\
\overline W&\geq A(1-e^{-\ell\mu})\geq M_{0}
\quad\text{on }(\overline\Omega\setminus B_{r}(x_{0}))\times[-T,0],\\
\overline W&\geq\frac{M_{0}}{h^{2}}(t-t_{0})^{2}\geq M_{0}
\quad\text{for }(x,t)\in\overline{\Omega_{T}}\text{ with }|t-t_{0}|\geq h.
\end{align*}
These inequalities give $\overline W\geq\varphi+\epsilon$ on
$\partial_{p}\Omega_{T}$, and $\overline W(z_{0})=m+\sigma$.
Define
\begin{align*}
\overline V:=\min\{M_{0},\overline W\}\quad\text{on }\overline{\Omega_{T}}.
\end{align*}
Since $M_{0}$ and $\overline W$ are viscosity supersolutions in $\Omega_{T}$, so is $\overline V$.
$\overline V\in C(\overline{\Omega_{T}})$,
and the estimates in \eqref{super es} follow from
the definition of $\overline V$ and the preceding bounds.

\smallskip
\noindent\textit{Case II. The initial boundary.}
Suppose that $t_{0}=-T$ and $x_{0}\in\overline\Omega$.
By continuity of $\varphi$ at $z_{0}$, choose $r>0$ and $0<\tau\leq T$ sufficiently small that
\begin{align*}
|\varphi(z)-\varphi(z_{0})|<\sigma/2
\quad\text{for }z\in\partial_{p}\Omega_{T}
\cap\bigl(\overline{B_{r}(x_{0})}\times[-T,-T+\tau]\bigr).
\end{align*}
We consider
\begin{align*}
\overline W(x,t):=m+\sigma+\frac{M_{0}}{r^{2}}|x-x_{0}|^{2}
+B(t+T)\quad\text{for }(x,t)\in\overline{\Omega_{T}},
\end{align*}
where $B>0$ is to be determined.
Since $|D\ov{W}|=2M_{0}|x-x_{0}|/r^{2}$, $D^{2}\overline W=2M_{0}I/r^{2}$, and $0<\nu<1$,
it suffices to choose $B$ so that
\begin{align*}
B\geq\frac{2n\Lambda M_{0}}{r^{2}}
\left(1+\frac{4M_{0}^{2}(\operatorname{diam}\Omega)^{2}}{r^{4}}\right)^{\kappa/2},\qquad B\tau\geq M_{0}.
\end{align*}
Indeed, by \ref{a1}, in $\Omega_{T}$ we have
\begin{align*}
a_{\nu}(D\overline W)F_{\zeta}(D^{2}\overline W)
-(\overline W)_{t}
&\leq a_{\nu}(D\overline W)\M^{+}(D^{2}\overline W)
-(\overline W)_{t} \\
&\leq \frac{2n\Lambda M_{0}}{r^{2}}
\left(1+\frac{4M_{0}^{2}(\operatorname{diam}\Omega)^{2}}{r^{4}}\right)^{\kappa/2}-B \\
&\leq0\leq g_{\epsilon}(\overline W).
\end{align*}
Thus $\overline W$ is a viscosity supersolution in $\Omega_{T}$.
Also,
\begin{align*}
\ov{W}&\geq m+\sigma/2\geq \varphi+\epsilon
\quad\text{on } \partial_{p}\Omega_{T}
\cap\bigl(\overline{B_{r}(x_{0})}\times[-T,-T+\tau]\bigr),\\
\overline W&\geq \frac{M_{0}}{r^{2}}|x-x_{0}|^{2}\geq M_{0}
\quad\text{on }(\overline\Omega\setminus B_{r}(x_{0}))\times[-T,0],\\
\overline W&\geq B(t+T)\geq B\tau\geq M_{0}
\quad\text{for }(x,t)\in\overline{\Omega_{T}}\text{ with }t\ge-T+\tau.
\end{align*}
Hence $\overline W\geq\varphi+\epsilon$ on $\partial_{p}\Omega_{T}$,
and $\overline W(z_{0})=m+\sigma$.
Define
\begin{align*}
\overline V:=\min\{M_{0},\overline W\}\quad\text{on }\overline{\Omega_{T}}.
\end{align*}
As in Case I, $\ov{V}\in C(\ov{\Omega_{T}})$ is a viscosity supersolution with \eqref{super es}.

In both Case I and Case II, $\overline W$ has a modulus of continuity near $z_{0}$ independent of $\epsilon,\nu,\zeta$, and so does $\overline V$.
\end{proof}

We next construct the corresponding subsolution barriers
using power profiles with exponent $\alpha$.

\begin{lem}
\label{bdrysub}
Under the assumptions of Lemma~\ref{bdrysuper}, for every
$z_{0}\in\partial_{p}\Omega_{T}$, $\sigma>0$, and
$\epsilon,\nu,\zeta\in(0,1)$, there exists a viscosity subsolution
$\underline V\in C(\overline{\Omega_{T}})$ to 
\eqref{pe} such that
\begin{align}
\label{subsol es}
\epsilon&\leq\underline V\leq M_{0}
\quad\text{in }\overline{\Omega_{T}},\qquad
\underline V\leq\varphi+\epsilon
\quad\text{on }\partial_{p}\Omega_{T}, \qquad \varphi(z_{0})+\epsilon-2\sigma\leq\underline V(z_{0})\leq\varphi(z_{0})+\epsilon.
\end{align}
Moreover, $\underline V$ has a modulus of continuity independent of $\epsilon,\nu,\zeta$ in a fixed relative neighborhood of $z_{0}$.
\end{lem}

\begin{proof}
Fix $z_{0}=(x_{0},t_{0})\in\partial_{p}\Omega_{T}$ and $\sigma>0$,
and write $m:=\varphi(z_{0})+\epsilon$.
If $m\leq2\sigma$, take $\underline V\equiv\epsilon$.
Since $g_{\epsilon}(\epsilon)=0$ and
$m-2\sigma\leq\epsilon\leq m$, all the assertions hold.
We therefore assume that $m>2\sigma$.
Let $\alpha$ and $\phi$ be as in Lemma~\ref{barrier homo}.
We use
\begin{align}
\label{basic ineq}
1<\alpha<2,\qquad
\kappa(\alpha-1)+\alpha-2=\alpha(\gamma-1),\qquad
\alpha-1-\alpha(\gamma-1)>0.
\end{align}
We divide the proof into two cases according to the location of $z_{0}$.

\smallskip
\noindent\textit{Case I. The lateral boundary.}
Suppose that $x_{0}\in\partial\Omega$ and $-T<t_{0}\leq0$.
Choose $y,\rho,\eta,r,h,t_{+}$, and $\mu$ as in Case I of
the proof of Lemma~\ref{bdrysuper}, for the same $z_{0}$ and $\sigma$.
Define
\begin{align*}
H(t)&:=(m-\sigma)^{1/\alpha}
\left(1-\frac{(t-t_{0})^{2}}{h^{2}}\right)
\quad\text{for }t\in[-T,0],\\
q(x,t)&:=H(t)-L\eta(x),\qquad
\underline W(x,t):=[q(x,t)]_{+}^{\alpha}
\quad\text{for }(x,t)\in\overline{\Omega_{T}},
\end{align*}
where $L>0$ is to be determined.
For $t\in[-T,0]$ with $H(t)>0$, we have
\begin{align}
\label{H bound}
0<H(t)\leq M_{0}^{1/\alpha},\qquad
|H'(t)|\leq\frac{2M_{0}^{1/\alpha}}{h}.
\end{align}
Choose $L$ so that
\begin{align*}
L\mu\geq M_{0}^{1/\alpha},\quad
L\geq\frac{2\Lambda(n-1)M_{0}^{1/\alpha}}
{\lambda(\alpha-1)\rho},\quad
\frac{\alpha^{\kappa+1}\lambda(\alpha-1)}{2}L^{\kappa+2}
\geq1+\frac{2\alpha M_{0}^{2-\gamma}}{h}.
\end{align*}
Fix $(x,t)\in\Omega_{T}$ with $q(x,t)>0$. Then
\begin{align*}
F_{\zeta}(D^{2}\underline W)
\geq \M^{-}(D^{2}\underline W)
\geq\alpha q^{\alpha-2}
\left(\lambda(\alpha-1)L^{2}
-\frac{\Lambda(n-1)Lq}{|x-y|}\right)
\geq\frac{\alpha\lambda(\alpha-1)}{2}L^{2}q^{\alpha-2}>0.
\end{align*}
Here we used \ref{a1}, Lemma~\ref{Hessian}, $q\leq M_{0}^{1/\alpha}$, $|x-y|\geq \rho$, and the choice of $L$.
From $a_{\nu}(D\underline W)\geq|D\underline W|^{\kappa}$, $|D\underline{W}|=\alpha Lq^{\alpha-1}$, and \eqref{basic ineq}, we obtain
\begin{align*}
a_{\nu}(D\underline W)F_{\zeta}(D^{2}\underline W)
\geq\frac{\alpha^{\kappa+1}\lambda(\alpha-1)}{2}
L^{\kappa+2}q^{\alpha(\gamma-1)}.
\end{align*}
Moreover, using \eqref{H bound}, $q\leq M_{0}^{1/\alpha}$, and \eqref{basic ineq}, we have
\begin{align*}
(\underline W)_{t}=\alpha q^{\alpha-1}H'
\leq \alpha M_{0}^{\frac{\alpha(2-\gamma)-1}{\alpha}}q^{\alpha(\gamma-1)}|H'| 
\leq\frac{2\alpha M_{0}^{2-\gamma}}{h}q^{\alpha(\gamma-1)}.    
\end{align*}
Hence, by the choice of $L$, $\underline{W}$ satisfies
\begin{align*}
a_{\nu}(D\underline{W})F_{\zeta}(D^{2}\underline{W})-(\underline{W})_{t}\geq q^{\alpha(\gamma-1)}=\underline{W}^{\gamma-1}\geq g_{\epsilon}(\underline{W}) \quad \text{in } \{\underline W>0\}\cap\Omega_{T}.    
\end{align*}
It remains to verify the subsolution property on $\{\underline W=0\}\cap\Omega_{T}$. 
Let $\xi\in C^{2,1}$ touch $\underline W$ from above at $z=(x,t)\in\Omega_{T}$ with $\xi(z)=\underline W(z)=0$.
Since $\xi\geq\underline W\geq0$ in a neighborhood of $z$, we have
\begin{align}
\label{min point}
D\xi(z)=0,\qquad D^{2}\xi(z)\geq0,\qquad \xi_{t}(z)\leq0,
\end{align}
where the last inequality follows by considering times less than $t$, which also covers the case $t=0$.
By \ref{a1}, \eqref{min point}, and $g_{\epsilon}(0)=0$, we obtain
\begin{align*}
a_{\nu}(D\xi(z))F_{\zeta}(D^{2}\xi(z))-\xi_{t}(z)
\geq0=g_{\epsilon}(\underline W(z)).
\end{align*}
Therefore $\underline W$ is a viscosity subsolution in $\Omega_{T}$.
By construction, we have
\begin{align*}
\underline W&\leq m-\sigma\leq\varphi+\epsilon
\quad\text{on }\partial_{p}\Omega_{T}
\cap\bigl(\overline{B_{r}(x_{0})}\times[t_{0}-h,t_{+}]\bigr),\\
\underline W&=0
\quad\text{on }(\overline\Omega\setminus B_{r}(x_{0}))\times[-T,0],\\
\underline W&=0
\quad\text{for }(x,t)\in\overline{\Omega_{T}}\text{ with }|t-t_{0}|\geq h.
\end{align*}
The second line follows from \eqref{eta bound} and
the choice of $L$, since
\begin{align*}
q(x,t)
\leq (m-\sigma)^{1/\alpha}-L\mu
\leq M_{0}^{1/\alpha}-L\mu
\leq 0 \qquad \text{for }(x,t)\in(\overline\Omega\setminus B_{r}(x_{0}))\times[-T,0].
\end{align*}
Thus $\underline W\leq\varphi+\epsilon$ on $\partial_{p}\Omega_{T}$,
and $\underline W(z_{0})=m-\sigma$.
Define
\begin{align*}
\underline V:=\max\{\epsilon,\underline W\}\quad\text{on }\overline{\Omega_{T}}.
\end{align*}
Since $\epsilon$ and $\underline W$ are viscosity subsolutions in $\Omega_{T}$, so is $\underline V$.
$\underline V\in C(\overline{\Omega_{T}})$ and satisfies \eqref{subsol es}.

\smallskip
\noindent\textit{Case II. The initial boundary.}
Suppose that $t_{0}=-T$ and $x_{0}\in\overline\Omega$.
Extend $\phi$ by zero to $\R^{n}$, so that $\phi\in C_{c}^{1}(\R^{n})$.
For $r>0$, consider
\begin{align*}
\underline W(x,t):=
\begin{cases}
(m-\sigma)\sqrt{1-B(t+T)}\,
\phi\!\left(\dfrac{x-x_{0}}{r\sqrt{1-B(t+T)}}\right),
& B(t+T)<1,\\
0,& B(t+T)\geq1,
\end{cases}
\end{align*}
for $(x,t)\in\overline{\Omega_{T}}$, where $B>0$ is to be determined.
$\underline W$ is continuous on $\overline{\Omega_{T}}$,
including at points where $B(t+T)=1$, since
\begin{align*}
0\leq\underline W
\leq(m-\sigma)\sqrt{[1-B(t+T)]_{+}}.
\end{align*}
By continuity of $\varphi$ at $z_{0}$, choose $r>0$ and $0<\tau\leq T$ small enough that
\begin{align*}
|\varphi(z)-\varphi(z_{0})|<\sigma/2
\quad\text{for }z\in\partial_{p}\Omega_{T}
\cap\bigl(\overline{B_{r}(x_{0})}\times[-T,-T+\tau]\bigr).
\end{align*}
Let $K_{0}$ and $k_{0}$ be the constants from Lemma \ref{barrier homo}. 
Decrease $r$, if necessary, so that
\begin{align*}
r\leq\sigma,\qquad
k_{0}r^{-(\kappa+2)}\sigma^{\kappa+2-\gamma}\geq1,
\end{align*}
and choose
\begin{align*}
B\frac{r^{\kappa+2}}{2M_{0}^{\kappa}}\geq K_{0},\qquad B\tau\geq1.
\end{align*}
Fix $(x,t)\in\{\underline W>0\}\cap\Omega_{T}$ and set
\begin{align*}
y:=\frac{x-x_{0}}{r\sqrt{1-B(t+T)}}\in B_{1}.
\end{align*}
Throughout the following calculations, $\phi$ and its derivatives
are evaluated at $y$.
Direct differentiation gives
\begin{align*}
D\underline W&=\frac{m-\sigma}{r}D\phi,\qquad
D^{2}\underline W=\frac{m-\sigma}{r^{2}\sqrt{1-B(t+T)}}D^{2}\phi,\\
-(\underline W)_{t}
&=\frac{B(m-\sigma)}{2\sqrt{1-B(t+T)}}\bigl(\phi-y\cdot D\phi\bigr)
\geq\frac{B(m-\sigma)}{2\sqrt{1-B(t+T)}}\phi,
\end{align*}
where we used $y\cdot D\phi(y)=-2\alpha|y|^{2}(1-|y|^{2})^{\alpha-1}\leq0$.
By the definition of $a_{\nu}$ and \ref{a1},
\begin{align*}
a_{\nu}(D\underline W)&=\left(\frac{m-\sigma}{r}\right)^{\kappa}a_{\nu r/(m-\sigma)}(D\phi),\\
F_{\zeta}(D^{2}\underline W)&\geq\M^{-}(D^{2}\underline W)
=\frac{m-\sigma}{r^{2}\sqrt{1-B(t+T)}}\M^{-}(D^{2}\phi).
\end{align*}
The choice of $B$ and $m-\sigma\leq M_{0}$ give
$$\frac{Br^{\kappa+2}}{2(m-\sigma)^{\kappa}}\geq K_{0}.$$
Since $r\leq\sigma<m-\sigma$ and $0<\nu<1$, we have $0<\nu r/(m-\sigma)<1$.
Combining these estimates and applying Lemma~\ref{barrier homo} with $\nu$ replaced by $\nu r/(m-\sigma)$, we obtain
\begin{align*}
a_{\nu}(D\underline W)F_{\zeta}(D^{2}\underline W)-(\underline W)_{t}
&\geq\frac{(m-\sigma)^{\kappa+1}}{r^{\kappa+2}\sqrt{1-B(t+T)}}
\left(a_{\nu r/(m-\sigma)}(D\phi)\M^{-}(D^{2}\phi)+\frac{Br^{\kappa+2}}{2(m-\sigma)^{\kappa}}\phi\right)\\
&\geq\frac{(m-\sigma)^{\kappa+1}}{r^{\kappa+2}\sqrt{1-B(t+T)}}
\bigl(a_{\nu r/(m-\sigma)}(D\phi)\M^{-}(D^{2}\phi)+K_{0}\phi\bigr)\\
&\geq\frac{k_{0}(m-\sigma)^{\kappa+1}}{r^{\kappa+2}\sqrt{1-B(t+T)}}\phi^{\gamma-1}
=\frac{k_{0}(m-\sigma)^{\kappa+2-\gamma}}{r^{\kappa+2}(1-B(t+T))^{\gamma/2}}\underline W^{\gamma-1}\\
&\geq\underline W^{\gamma-1}\geq g_{\epsilon}(\underline W).
\end{align*}
For the last line, we used $0<1-B(t+T)\leq1$, $m-\sigma>\sigma$, and the choice of $r$, which give
\begin{align*}
\frac{k_{0}(m-\sigma)^{\kappa+2-\gamma}}{r^{\kappa+2}(1-B(t+T))^{\gamma/2}}
\geq k_{0}r^{-(\kappa+2)}\sigma^{\kappa+2-\gamma}\geq1,
\end{align*}
together with $g_{\epsilon}(s)\leq s^{\gamma-1}$ for every $s>0$.
Thus $\underline W$ satisfies the subsolution inequality in $\{\underline W>0\}\cap\Omega_{T}$.
On $\left\{\underline W=0\right\}\cap \Omega_{T}$, the same argument as in Case I applies.
Hence $\underline W$ is a viscosity subsolution in $\Omega_{T}$.
Moreover, we have
\begin{align*}
\underline W&\leq m-\sigma\leq\varphi+\epsilon
\quad\text{on }\partial_{p}\Omega_{T}
\cap\bigl(\overline{B_{r}(x_{0})}\times[-T,-T+\tau]\bigr),\\
\underline W&=0
\quad\text{on }(\overline\Omega\setminus B_{r}(x_{0}))\times[-T,0],\\
\underline W&=0
\quad\text{for }(x,t)\in\overline{\Omega_{T}}\text{ with } t\geq -T+\tau.
\end{align*}
Hence $\underline W\leq\varphi+\epsilon$ on $\partial_{p}\Omega_{T}$,
and $\underline W(z_{0})=m-\sigma$.
Define
\begin{align*}
\underline V:=\max\{\epsilon,\underline W\}\quad\text{on }\overline{\Omega_{T}}.
\end{align*}
As in Case I, $\underline V$ is a continuous viscosity subsolution in $\Omega_{T}$ satisfying \eqref{subsol es}.

Finally, note that $m-\sigma\leq M_{0}$ and that the auxiliary parameters in both constructions can be chosen independently of $\epsilon,\nu,\zeta$.
Therefore, $\underline V$ has a modulus of continuity independent of $\epsilon,\nu,\zeta$ in a fixed relative neighborhood of $z_{0}$.
This proves the lemma.
\end{proof}

Comparing any solution of \eqref{pe} with these barriers
yields the following boundary continuity estimate,
uniform in $\epsilon$, $\nu$, and $\zeta$.

\begin{cor}
\label{bdrycont}
Under the assumptions of Lemma~\ref{bdrysuper}, for each
$\epsilon,\nu,\zeta\in(0,1)$, let
$u_{\epsilon,\nu,\zeta}\in C(\overline{\Omega_{T}})$ be a viscosity
solution to \eqref{pe}.
Then, for every $z_{0}\in\partial_{p}\Omega_{T}$ and $\omega>0$,
there exists a relative neighborhood $N$ of $z_{0}$ in
$\overline{\Omega_{T}}$, independent of $\epsilon,\nu,\zeta$,
such that
\begin{align*}
|u_{\epsilon,\nu,\zeta}(z)-\varphi(z_{0})-\epsilon|<\omega
\qquad\text{for }z\in N\cap\Omega_{T}.
\end{align*}
\end{cor}

\begin{proof}
Fix $z_{0}\in\partial_{p}\Omega_{T}$ and $\omega>0$,
and choose $0<\sigma<\omega/4$.
Let $\overline V$ and $\underline V$ be the barriers constructed
in Lemmas~\ref{bdrysuper} and \ref{bdrysub}, respectively,
for this choice of $z_{0}$ and $\sigma$.
By Lemma~\ref{comparison} and continuity at $t=0$, we have
\begin{align*}
\underline V\leq u_{\epsilon,\nu,\zeta}\leq\overline V
\quad\text{in }\Omega_{T}.
\end{align*}
The uniform moduli of continuity of $\underline V$ and
$\overline V$ yield a relative neighborhood $N$ of $z_{0}$ in
$\overline{\Omega_{T}}$, independent of $\epsilon,\nu,\zeta$,
such that
\begin{align*}
|\underline V(z)-\underline V(z_{0})|<\sigma,\qquad
|\overline V(z)-\overline V(z_{0})|<\sigma
\quad\text{for }z\in N\cap\Omega_{T}.
\end{align*}
Consequently, by \eqref{super es} and \eqref{subsol es},
\begin{align*}
u_{\epsilon,\nu,\zeta}(z)
&\geq\underline V(z)>\underline V(z_{0})-\sigma
\geq\varphi(z_{0})+\epsilon-3\sigma,\\
u_{\epsilon,\nu,\zeta}(z)
&\leq\overline V(z)<\overline V(z_{0})+\sigma
\leq\varphi(z_{0})+\epsilon+2\sigma
\end{align*}
for $z\in N\cap\Omega_{T}$. Hence
\begin{align*}
|u_{\epsilon,\nu,\zeta}(z)-\varphi(z_{0})-\epsilon|
<3\sigma<\omega,
\end{align*}
which completes the proof.
\end{proof}

\subsection{Existence and interior smoothness}
The comparison principle and boundary barriers allow us
to apply Perron's method to \eqref{pe}.
For fixed parameters, we apply the gradient regularity
result and the Schauder estimate from Section~\ref{sec2}
and then use a bootstrap argument to obtain the interior
smoothness needed for Bernstein's method.

\begin{thm}
\label{existence}
Let $\kappa>0$, $0<\gamma<1$, and assume \ref{a1}--\ref{a3}.
Let $\varphi\in C(\partial_{p}\Omega_{T})$ be nonnegative.
For every $\epsilon,\nu,\zeta\in(0,1)$, there exists a unique
viscosity solution
$u_{\epsilon,\nu,\zeta}\in C(\overline{\Omega_{T}})
\cap C^{\infty}(\Omega_{T})$ to \eqref{pe} satisfying
\begin{align}
\label{es}
\epsilon\leq u_{\epsilon,\nu,\zeta}
\leq\norm{\varphi}_{L^{\infty}(\partial_{p}\Omega_{T})}+\epsilon
\quad\text{in }\overline{\Omega_{T}}.
\end{align}
\end{thm}

\begin{proof}
Fix $\epsilon,\nu,\zeta\in(0,1)$ and set
$M=\norm{\varphi}_{L^{\infty}(\partial_{p}\Omega_{T})}$.
To treat $t=0$ as an interior time, extend the lateral boundary
data by setting $\varphi(x,t)=\varphi(x,0)$ for
$(x,t)\in\partial\Omega\times(0,1)$, keeping the same notation.
Let
\begin{align*}
E:=\Omega\times(-T,1), \qquad \Gamma:=(\overline\Omega\times\{-T\})
\cup(\partial\Omega\times[-T,1)).
\end{align*}
For $z_0\in\Gamma$ and $\sigma>0$, let $\overline V_{z_0,\sigma},\underline V_{z_0,\sigma} \in C(E\cup\Gamma)$ denote the supersolution and subsolution in $E$, respectively, constructed in Lemmas~\ref{bdrysuper} and \ref{bdrysub} for the extended boundary data.
Applying Lemma~\ref{comparison} with
$H(q,X)=a_{\nu}(q)F_{\zeta}(X)$ and $f=g_{\epsilon}$, we obtain
\begin{align}
\label{perron comp}
\underline V_{z_1,\sigma_1}
\leq \overline V_{z_2,\sigma_2}
\quad\text{in }E\cup\Gamma
\end{align}
for all $z_1,z_2\in\Gamma$ and $\sigma_1,\sigma_2>0$.
Set
\begin{align*}
\underline v
&:=\bigl(\sup\nolimits_{z_0,\sigma}
\underline V_{z_0,\sigma}\bigr)^*,
\qquad
\overline v
:=\bigl(\inf\nolimits_{z_0,\sigma}
\overline V_{z_0,\sigma}\bigr)_*,
\end{align*}
where the semicontinuous envelopes are taken relative to
$E\cup\Gamma$.
By stability, $\underline v$ and $\overline v$ are respectively a viscosity subsolution and a viscosity supersolution in $E$.
For fixed $z_0\in\Gamma$ and $\sigma>0$, \eqref{perron comp} and the continuity of $\overline V_{z_0,\sigma}$ and $\underline V_{z_0,\sigma}$ give
\begin{align*}
\underline V_{z_0,\sigma}
\leq w\leq
\overline V_{z_0,\sigma}
\quad\text{in }E\cup\Gamma,
\qquad
w\in\{\underline v,\overline v\}.
\end{align*}
Letting $z\to z_0$ and then $\sigma\to0^+$, we obtain from
\eqref{super es} and \eqref{subsol es} that
\begin{align*}
\lim_{\substack{z\to z_0\\z\in E\cup\Gamma}}\underline v(z)
=
\lim_{\substack{z\to z_0\\z\in E\cup\Gamma}}\overline v(z)
=
\varphi(z_0)+\epsilon
\end{align*}
for every $z_0\in\Gamma$.
Lemma~\ref{comparison} then yields
$\underline v\leq\overline v$ in $E$.
Perron's method (see \cite[Theorem 4.1 and Section 8]{CIL1992}),
with comparison provided by Lemma~\ref{comparison},
yields a viscosity solution $u\in C(E\cup\Gamma)$ to
\begin{align}
\label{extended}
\left\{
\begin{aligned}
a_{\nu}(Du)F_{\zeta}(D^{2}u)-u_t
&=g_{\epsilon}(u)
&&\text{in }E,\\
u&=\varphi+\epsilon
&&\text{on }\Gamma.
\end{aligned}
\right.
\end{align}
Since $u\geq\underline v\geq\epsilon$, comparison with the constant supersolution $M+\epsilon$ gives $\epsilon\leq u\leq M+\epsilon$ in $E\cup\Gamma$.
The restriction $u_{\epsilon,\nu,\zeta}:=u|_{\overline{\Omega_{T}}}$
solves \eqref{pe} and satisfies \eqref{es}.
Uniqueness follows from Lemma~\ref{comparison} and continuity
at $t=0$.

We next prove that $u\in C^{\infty}(E)$.
Define
\begin{align*}
U(x,z,t)&:=u(x,t)+\nu z \quad \text{for } (x,z,t)\in\Omega\times(-1,1)\times(-T,1), \\ \widetilde F_{\zeta}(N)
&:=F_{\zeta}(N_{xx})+\lambda N_{zz}
\quad \text{for } N\in\mathcal S^{n+1},
\end{align*}
where $N_{xx}$ denotes the upper-left $n\times n$ block of $N$.
$\widetilde F_{\zeta}$ satisfies \ref{a1}
and \ref{a2} in dimension $n+1$ with the same ellipticity
constants. $\widetilde F_{\zeta}$ is smooth
and $D\widetilde F_{\zeta}$ is globally Lipschitz continuous.
We claim that $U$ is a viscosity solution to
\begin{align}
\label{lifted-equation}
|D_{x,z}U|^{\kappa}
\widetilde F_{\zeta}(D^2_{x,z}U)
-U_{t}=g_{\epsilon}(u(x,t)) \quad \text{in }\Omega\times(-1,1)\times(-T,1).
\end{align}
Let $\Psi\in C^{2,1}$ touch $U$ from above at
$(x_0,z_0,t_0)\in \Omega\times(-1,1)\times(-T,1)$. Then $\Psi_z=U_{z}=\nu$ and $\Psi_{zz}\geq0$ at $(x_{0}, z_{0}, t_{0})$. Since $\Psi(x,z_0,t)-\nu z_0$ touches $u$ from above at $(x_0,t_0)$ and $u$ is a viscosity solution to \eqref{extended}, we have
\begin{align*}
g_{\epsilon}(u(x_{0}, t_{0}))&\leq a_{\nu}(D_{x}\Psi)F_{\zeta}(D_{x}^{2}\Psi)-\Psi_{t}\\
&=(\nu^{2}+|D_{x}\Psi|^{2})^{\kappa/2}F_{\zeta}(D_{x}^{2}\Psi)-\Psi_{t} \\
&\leq |D_{x,z}\Psi|^{\kappa}\widetilde{F}_{\zeta}(D^{2}_{x,z}\Psi)-\Psi_{t},
\end{align*}
where all derivatives of $\Psi$ are evaluated at $(x_0,z_0,t_0)$, and the last inequality follows from $\Psi_z=\nu$ and $\widetilde{F}_{\zeta}(D^{2}_{x,z}\Psi) =F_{\zeta}(D_x^2\Psi)+\lambda\Psi_{zz} \geq F_{\zeta}(D_x^2\Psi)$.

This proves that $U$ is a viscosity subsolution to \eqref{lifted-equation}. The supersolution property follows analogously by considering test functions that touch $U$ from below, thereby proving the claim.
Since $g_{\epsilon}(u)$ is bounded and continuous,
Lemma~\ref{gradient-regularity}, applied in dimension $n+1$,
gives $U\in C^{1+\sigma,(1+\sigma)/2}_{\loc}(\Omega\times(-1,1)\times(-T,1))$, which implies
\begin{align*}
u\in C^{1+\sigma,(1+\sigma)/2}_{\loc}(E), Du\in C_{\loc}^{\sigma, \sigma/2}(E),
\end{align*}
for some $\sigma\in(0,1)$.
Thus $a_{\nu}(Du)$ and $g_{\epsilon}(u)$ are locally Hölder continuous. Fix a cylinder $Q\Subset E$. Then, in $Q$,
\begin{align*}
0<\nu^{\kappa}
\leq a_{\nu}(Du)
\leq
\left(\nu^2+\norm{Du}_{L^{\infty}(Q)}^2\right)^{\kappa/2}.
\end{align*}
Consequently, $G(x,t,X):=a_{\nu}(Du(x,t))F_{\zeta}(X)$ is uniformly elliptic in $Q$, is convex or concave in $X$, and satisfies \eqref{coeff holder}.
By Lemma~\ref{schauder} with $f(x,t)=-g_{\epsilon}(u(x,t))$, we have
\begin{align*}
u\in C^{2+\eta,1+\eta/2}_{\loc}(Q)
\end{align*}
for some $\eta\in(0,1)$. 
Since $a_{\nu}$, $F_{\zeta}$, and $g_{\epsilon}$ are smooth, applying the linear interior Schauder estimates to difference quotients in space and time gives $u_{x_k},u_t\in C^{2+\eta,1+\eta/2}_{\loc}(Q)$.
A bootstrap argument then yields $u\in C^{\infty}(Q)$. 
Since $Q\Subset E$ was arbitrary, $u\in C^{\infty}(E)$. Together with the continuity obtained above, this yields $u_{\epsilon,\nu,\zeta}\in C(\overline{\Omega_{T}})\cap C^{\infty}(\Omega_{T})$.
\end{proof}

\section{Local compactness}
\label{sec4}
We establish the local compactness needed to pass to the
limit in \eqref{pe} as $(\nu,\zeta)\to(0,0)$ with $\epsilon$
fixed, and then in \eqref{pe1} as $\epsilon\to0$.
Together with the uniform boundary estimates, this yields
a viscosity solution to \eqref{eq}.

\subsection{A weak Bernstein estimate}
We first estimate $|Du|^{\kappa+2}/u^a$ for a suitable
$a\in(0,\gamma)$ by Bernstein's method.
The estimate is uniform in the regularization parameters
$\epsilon$, $\nu$, and $\zeta$.

\begin{lem}\label{bernstein}
Let $\kappa>0$ and $0<\gamma<1$, and set
\begin{align*}
 C_{0}:=\Lambda+\frac{(\kappa+2)n\Lambda^{2}}{\lambda},
 \qquad
 a:=\min\left\{\frac{\gamma}{2},\frac{\lambda}{4C_{0}}\right\}.
\end{align*}
Suppose that $F\in C^{\infty}(\mathcal{S}^{n})$ satisfies \ref{a1}.
Let $0\leq\beta_{\epsilon}\leq1$ be smooth and nondecreasing.
If $0<u\leq M$ is a smooth solution to
\begin{align*}
(\nu^{2}+|Du|^{2})^{\kappa/2}F(D^{2}u)-u_{t}
 =\delta\beta_{\epsilon}(u)u^{\gamma-1}
 \quad\text{in }B_{r}\times(-T,0],
\end{align*}
where $\nu>0$ and $\delta\geq0$, then
\begin{align}\label{weak-bernstein-estimate}
 \sup_{B_{r/2}\times[-T/2,0]}\frac{|Du|^{\kappa+2}}{u^{a}}
 \leq C\left(
 \delta M^{\gamma-a}
 +\frac{M^{\kappa+2-a}}{r^{\kappa+2}}
 +\frac{M^{2-a}}{T}\right),
\end{align}
where $C=C(n,\lambda,\Lambda, \kappa,\gamma)>0$.
\end{lem}

\begin{proof}
Note that \ref{a1} implies
\begin{align}
\label{diff uniform ellip}
 \lambda I\leq DF(X)\leq\Lambda I \quad \text{for every } X\in\mathcal{S}^{n}.
\end{align}
Let $b:=a/(\kappa+2)$ and write
\begin{align*}
 \rho(x,t):=(\nu^{2}+|Du(x,t)|^{2})^{1/2},\qquad
 F_{ij}(x,t):=(DF(D^{2}u(x,t)))_{ij},\qquad
 g(s):=\delta\beta_{\epsilon}(s)s^{\gamma-1}.
\end{align*}
Then we have 
\begin{align}
\label{rewrite eq}
\rho^{\kappa}F(D^{2}u)-u_{t}=g(u).    
\end{align}
Define the linear operator $\mathcal{P}$ by
\begin{align*}
 \mathcal{P}v
 :=v_{t}-\rho^{\kappa}F_{ij}v_{ij}
 -\kappa\rho^{\kappa-2}F(D^{2}u)u_{i}v_{i}.
\end{align*}
For smooth functions $h$ and $v$, the chain rule gives
\begin{align}
\label{chain rule}
\mathcal{P}(h(v))&=h'(v)\mathcal{P}v-\rho^{\kappa}F_{ij}h''(v)v_{i}v_{j}.
\end{align}
Equation \eqref{rewrite eq}, its spatial derivatives, and the chain rule
\eqref{chain rule} give
\begin{align*}
 \mathcal{P}u_{k}&=-g'(u)u_{k},\\
 \mathcal{P}u
 &=\rho^{\kappa}\left[
 \left(1-\frac{\kappa|Du|^{2}}{\rho^{2}}\right)F(D^{2}u)
 -F_{ij}u_{ij}\right]-g(u), \\
\mathcal{P}(|Du|^{2})
&=-2g'(u)|Du|^{2}-2\rho^{\kappa}F_{ij}u_{ik}u_{jk}.    
\end{align*}
Choose $\eta\in C^{\infty}(\R^{n+1})$ such that
\begin{align*}
 &0\leq\eta\leq1,\quad
 \eta=1\quad\text{in }B_{r/2}\times[-T/2,0], \quad \supp\eta\cap\partial_{p}(B_{r}\times(-T,0])=\varnothing,\\
 &|D\eta|\leq C/r,\qquad |D^{2}\eta|\leq C/r^{2},
 \qquad |\eta_{t}|\leq C/T.
\end{align*}
Define
\begin{align*}
 V:=\eta\frac{|Du|}{u^{b}}
 \quad\text{in }B_{r}\times(-T,0].
\end{align*}
Extending $V$ by zero to $\partial_{p}(B_{r}\times(-T,0])$, we have
$V\in C(\overline{B_{r}}\times[-T,0])$.
If $\max_{\overline{B_{r}}\times[-T,0]}V=0$, the conclusion is immediate.
Otherwise, let $(x_{0},t_{0})\in B_{r}\times(-T,0]$ satisfy
\begin{align*}
 V(x_{0},t_{0})=\max_{\overline{B_{r}}\times[-T,0]}V>0.
\end{align*}
Since $V(x_{0},t_{0})>0$, $\log V$ is well defined
near $(x_{0},t_{0})$ and attains a local maximum there. Hence
$D(\log V)=0$ and $D^{2}(\log V)\leq0$ at this point.
Moreover, $(\log V)_{t}=0$ if $t_{0}<0$.
If $t_{0}=0$, the maximality of $V$ gives
$(\log V)_{t}(x_{0},0)\geq0$, where the time derivative is taken
from the left. 
Consequently, we have at $(x_{0},t_{0})$
\begin{align}\label{bernstein-maximum-relation}
 \eta>0,\quad |Du|>0,\quad
 \frac{u_{ik}u_{k}}{|Du|^{2}}
 =b\frac{u_{i}}{u}-\frac{\eta_{i}}{\eta},\quad
 0\leq\mathcal{P}(\log V),
\end{align}
where we used \eqref{diff uniform ellip} in the last inequality.
All subsequent calculations are evaluated at $(x_{0},t_{0})$.

Using
\begin{align*}
 \mathcal{P}(\log v)
 =\frac{\mathcal{P}v}{v}
 +\frac{\rho^{\kappa}}{v^{2}}F_{ij}v_{i}v_{j}
 \qquad(v>0)
\end{align*}
and \eqref{bernstein-maximum-relation}, we group the terms as follows:
\begin{align}\label{bernstein-log-identity}
 \mathcal{P}(\log V)
 &=\mathcal{P}(\log\eta)
 +\frac{1}{2}\mathcal{P}\bigl(\log(|Du|^{2})\bigr)
 -b\mathcal{P}(\log u)
 \nonumber\\
 &=\mathrm{I}+\mathrm{II}+\mathrm{III},
\end{align}
where
\begin{align*}
 \mathrm{I}:={}&-\frac{\rho^{\kappa}}{|Du|^{2}}F_{ij}u_{ik}u_{jk}\\
 &+\frac{b\rho^{\kappa}}{u}\left[
 F_{ij}u_{ij}
 +\left(\frac{\kappa|Du|^{2}}{\rho^{2}}-1\right)F(D^{2}u)
 \right]
 \\
 &-\frac{\kappa\rho^{\kappa-2}}{\eta}F(D^{2}u)u_{i}\eta_{i},\\
 \mathrm{II}:={}&-(b-2b^{2})\frac{\rho^{\kappa}}{u^{2}}F_{ij}u_{i}u_{j}
 -\frac{4b\rho^{\kappa}}{u\eta}F_{ij}u_{i}\eta_{j},\\
 \mathrm{III}:={}&-g'(u)+\frac{bg(u)}{u}
 +\frac{\eta_{t}}{\eta}
 -\frac{\rho^{\kappa}}{\eta}F_{ij}\eta_{ij}
 +\frac{3\rho^{\kappa}}{\eta^{2}}F_{ij}\eta_{i}\eta_{j}.
\end{align*}
We first estimate $\mathrm{I}$, which collects the terms involving
$D^{2}u$. By uniform ellipticity and $F(0)=0$,
\begin{align*}
 F_{ij}u_{ik}u_{jk}\geq\lambda|D^{2}u|^{2},
 \qquad
 \max\{|F(D^{2}u)|,|F_{ij}u_{ij}|\}
 \leq\sqrt{n}\Lambda|D^{2}u|.
\end{align*}
Since $\rho\geq|Du|$ and
$1+\left|\kappa|Du|^{2}/\rho^{2}-1\right|\leq\kappa+2$, we have
\begin{align*}
\mathrm{I}
&\leq-\frac{\lambda\rho^{\kappa}}{|Du|^{2}}|D^{2}u|^{2}
 +\sqrt{n}\Lambda\rho^{\kappa}
 \left(\frac{(\kappa+2)b}{u}
       +\frac{\kappa|u_{i}\eta_{i}|}{\eta\rho^{2}}\right)|D^{2}u|\\
 &\leq-\frac{\lambda\rho^{\kappa}}{|Du|^{2}}|D^{2}u|^{2}
 +\sqrt{n}\Lambda\rho^{\kappa}
 \left(\frac{(\kappa+2)b}{u}
       +\frac{\kappa|D\eta|}{\eta|Du|}\right)|D^{2}u|\\
 &\leq
 \frac{(\kappa+2)^{2}n\Lambda^{2}}{2\lambda}b^{2}
 \frac{\rho^{\kappa}|Du|^{2}}{u^{2}}
 +\frac{\kappa^{2}n\Lambda^{2}}{2\lambda}
 \frac{\rho^{\kappa}|D\eta|^{2}}{\eta^{2}}.
\end{align*}
Here we used $|u_{i}\eta_{i}|\leq|Du||D\eta|$,
$-AX^{2}+BX\leq B^{2}/(4A)$ for $A>0$, and
$(A+B)^{2}\leq2A^{2}+2B^{2}$.
Next, we estimate $\mathrm{II}$.
Using \eqref{diff uniform ellip} and Young's inequality
$AB\leq\tau A^{2}+B^{2}/(4\tau)$ for $\tau>0$, we obtain
\begin{align*}
 -(b-2b^{2})\frac{\rho^{\kappa}}{u^{2}}F_{ij}u_{i}u_{j}
 &\leq(-\lambda b+2\Lambda b^{2})
       \frac{\rho^{\kappa}|Du|^{2}}{u^{2}},\\
 \frac{4b\rho^{\kappa}}{u\eta}|F_{ij}u_{i}\eta_{j}|
 &\leq\frac{\lambda b}{4}
       \frac{\rho^{\kappa}|Du|^{2}}{u^{2}}
 +\frac{16b\Lambda^{2}}{\lambda}
       \frac{\rho^{\kappa}|D\eta|^{2}}{\eta^{2}}.
\end{align*}
We now collect the coefficients of $\rho^{\kappa}|Du|^{2}/u^{2}$.
Using $\kappa>0$ and recalling the definition of $C_{0}$, our choice of $a$ gives
\begin{align*}
 \frac{(\kappa+2)^{2}n\Lambda^{2}}{2\lambda}b^{2}-\lambda b+2\Lambda b^{2}+\frac{\lambda b}{4}
 \leq(\kappa+2)C_{0}b^{2}-\frac{3\lambda b}{4}
 =C_{0}ab-\frac{3\lambda b}{4}\leq-\frac{\lambda b}{2}.
\end{align*}
Finally, using $\beta_{\epsilon}'\geq0$, we estimate the source terms
in $\mathrm{III}$ as follows:
\begin{align*}
 -g'(u)+\frac{bg(u)}{u}
 &=\delta u^{\gamma-2}
 \left[(1-\gamma+b)\beta_{\epsilon}(u)
       -u\beta_{\epsilon}'(u)\right]
 \leq C\delta u^{\gamma-2}.
\end{align*}
Using uniform ellipticity, the bounds for $\eta$, and $0<\eta\leq1$,
we obtain
\begin{align*}
 \mathrm{III}
 \leq C\delta u^{\gamma-2}
 +\frac{C\rho^{\kappa}}{r^{2}\eta^{2}}
 +\frac{C}{T\eta}.
\end{align*}
Combining these estimates with \eqref{bernstein-log-identity},
we find
\begin{align*}
 0\leq-\frac{\lambda b}{2}
         \frac{\rho^{\kappa}|Du|^{2}}{u^{2}}
 +C\delta u^{\gamma-2}
 +\frac{C\rho^{\kappa}}{r^{2}\eta^{2}}
 +\frac{C}{T\eta}.
\end{align*}
Dividing by $\rho^{\kappa}|Du|^{2}/u^{2}$, we have
\begin{align*}
0\leq -\frac{\lambda b}{2}+C\frac{\delta u^{\gamma}}{\rho^{\kappa}|Du|^{2}}+C\frac{u^{2}}{r^{2}\eta^{2}|Du|^{2}}+C\frac{u^{2}}{T\eta\rho^{\kappa}|Du|^{2}}.
\end{align*}
Since $\rho^{\kappa}|Du|^{2}\geq|Du|^{\kappa+2}$ and
$|Du|^{\kappa+2}=V^{\kappa+2}u^{a}\eta^{-(\kappa+2)}$, we get
\begin{align*}
 \frac{\delta u^{\gamma}}{\rho^{\kappa}|Du|^{2}}\leq\frac{\delta M^{\gamma-a}}{V^{\kappa+2}}, \qquad
 \frac{u^{2}}{r^{2}\eta^{2}|Du|^{2}}
 \leq\frac{M^{2-2b}}{r^{2}V^{2}}, \qquad
 \frac{u^{2}}{T\eta\rho^{\kappa}|Du|^{2}}
 \leq\frac{M^{2-a}}{TV^{\kappa+2}}.
\end{align*}
Consequently,
\begin{align*}
 0\leq-\frac{\lambda b}{2}
 +C\frac{M^{2-2b}}{r^{2}V^{2}}
 +C\frac{\delta M^{\gamma-a}+M^{2-a}/T}{V^{\kappa+2}}.
\end{align*}
It follows that
\begin{align*}
 \max_{\overline{B_{r}}\times[-T,0]}V^{\kappa+2}
 \leq C\left(
 \delta M^{\gamma-a}
 +\frac{M^{\kappa+2-a}}{r^{\kappa+2}}
 +\frac{M^{2-a}}{T}\right).
\end{align*}
Since $\eta=1$ in $B_{r/2}\times[-T/2,0]$, this proves
\eqref{weak-bernstein-estimate}.
\end{proof}

\begin{rmk}\label{bernstein-uniformity}
The proof uses neither convexity nor homogeneity of $F$.
In particular, \eqref{weak-bernstein-estimate} applies to
$F_{\zeta}$ with the same constant
$C=C(n,\lambda,\Lambda,\kappa,\gamma)$.
\end{rmk}

\subsection{Time compactness}

With $\epsilon$ fixed, we first obtain a local H\"older
estimate in time, uniform in $\nu$ and $\zeta$.
This allows us to pass from $u_{\epsilon,\nu,\zeta}$ to $u_\epsilon$.

\begin{lem}
\label{regularized time compactness}
Let $\kappa>0$, $0<\gamma<1$, and $a>0$.
Suppose that $F$ satisfies \ref{a1},
and let $Q'\Subset Q$ be cylinders.
Let $u$ be a smooth solution to
\begin{align}
\label{time derivative eq}
(\nu^{2}+|Du|^{2})^{\kappa/2}F(D^{2}u)-u_{t}
=\delta g_{\epsilon}(u)
\quad\text{in }Q,
\end{align}
where $\epsilon,\nu\in(0,1)$ and $\delta\in[0,1]$.
Assume that, for some $M,K>0$,
\begin{align}
\label{uniform gradient quotient}
0<u\leq M \quad\text{in }Q,
\qquad
\sup_{Q}\frac{|Du|^{\kappa+2}}{u^{a}}\leq K.
\end{align}
Then
\begin{align}
\label{regularized spacetime modulus}
|u(x,t)-u(y,s)|
\leq C_{\epsilon}\bigl(|x-y|+|t-s|^{1/2}\bigr)
\qquad\text{for }(x,t),(y,s)\in Q',
\end{align}
where $C_{\epsilon}>0$ depends only on
$n,\lambda,\Lambda,\kappa,\gamma,a,M,K,\epsilon$, and $Q',Q$.
\end{lem}

\begin{proof}
Let $L:=(KM^{a})^{1/(\kappa+2)}$, so that $|Du|\leq L$ in $Q$.
Also, $0\leq\delta g_{\epsilon}(u)\leq\epsilon^{\gamma-1}$ in $Q$.
Choose $C>2n\Lambda(1+L^{2})^{\kappa/2}$ depending only on
$n,\Lambda,\kappa$, and $L$.
Fix $(x,s),(x,t)\in Q'$ with $s<t$, and choose $R>0$,
depending only on $Q',Q$, such that
$\overline{B_{R}(x)}\times[s,t]\subset Q$.
For $\eta>0$, define
\begin{align*}
A_{\eta}:=\frac{L^{2}}{4\eta}+\frac{M}{R^{2}}, \qquad
h^{\pm}(y,\tau):=u(x,s)\pm\left\{
\eta+A_{\eta}|y-x|^{2}
+\bigl(CA_{\eta}+\epsilon^{\gamma-1}\bigr)(\tau-s)
\right\}.
\end{align*}
We claim that $h^{-}\leq u\leq h^{+}$ in $\overline{B_{R}(x)}\times[s,t]$.
First, we have $h^{-}\leq u\leq h^{+}$ on $\partial_{p}(B_{R}(x)\times(s,t])$.
Indeed, at $\tau=s$, we have
\begin{align*}
h^{+}(y,s)=u(x,s)+\eta+A_{\eta}|y-x|^{2}\geq u(x,s)+\eta+L^{2}|y-x|^{2}/(4\eta)\geq u(x,s)+L|y-x|\geq u(y,s), \\
h^{-}(y,s)=u(x,s)-\eta-A_{\eta}|y-x|^{2}\leq u(x,s)-\eta-L^{2}|y-x|^{2}/(4\eta)\leq u(x,s)-L|y-x|\leq u(y,s).
\end{align*}
On $\partial B_{R}(x)\times[s,t]$, using $A_{\eta}R^{2}\geq M$, we have
\begin{align*}
h^{+}(y,\tau)\geq u(x,s)+M\geq u(y,\tau), \qquad h^{-}(y,\tau)\leq u(x,s)-M\leq u(y,\tau).
\end{align*}
Assume that $u-h^{+}$ has a positive maximum at $z_{0}$ in $\overline{B_{R}(x)}\times[s,t]$.
The boundary inequalities imply that $z_{0}\in B_{R}(x)\times(s,t]$.
Then $D^{2}u\leq D^{2}h^{+}=2A_{\eta}I$ and $u_{t}\geq h^{+}_{t}=CA_{\eta}+\epsilon^{\gamma-1}$ at $z_{0}$.
We have $u_{t}=h^{+}_{t}$ when $\tau<t$, while the inequality at $\tau=t$ follows from the left-sided time derivative.
By \eqref{time derivative eq} and \ref{a1}, we obtain at $z_{0}$
\begin{align*}
u_{t}= (\nu^{2}+|Du|^{2})^{\kappa/2}F(D^{2}u)-\delta g_{\epsilon}(u)&\leq (\nu^{2}+|Du|^{2})^{\kappa/2}\M^{+}(D^{2}u)-\delta g_{\epsilon}(u) \\
&\leq (\nu^{2}+|Du|^{2})^{\kappa/2}\M^{+}(D^{2}h^{+}) \\
&\leq 2n\Lambda(1+L^{2})^{\kappa/2}A_{\eta}<h^{+}_{t},
\end{align*}
which is a contradiction.
Similarly, suppose that $h^{-}-u$ attains a positive maximum at $z_{1}$ in $\overline{B_{R}(x)}\times[s,t]$.
Then $D^{2}u\geq D^{2}h^{-}=-2A_{\eta}I$ and $u_{t}\leq h^{-}_{t}=-(CA_{\eta}+\epsilon^{\gamma-1})$ at $z_{1}$, whereas
\begin{align*}
u_{t}
\geq-2n\Lambda(1+L^{2})^{\kappa/2}A_{\eta}
-\epsilon^{\gamma-1}>h^{-}_{t},
\end{align*}
a contradiction.
Hence $h^{-}\leq u\leq h^{+}$ in
$\overline{B_{R}(x)}\times[s,t]$, and
\begin{align*}
|u(x,t)-u(x,s)|
\leq\eta+
\left[\frac{CL^{2}}{4\eta}+\frac{CM}{R^{2}}
+\epsilon^{\gamma-1}\right](t-s).
\end{align*}
For $0<t-s\leq1$, taking $\eta=\sqrt{t-s}$ gives $|u(x,t)-u(x,s)|\leq C_{\epsilon}|t-s|^{1/2}$.
For $t-s>1$, the same estimate follows directly from $0<u\leq M$.
Combining this with the spatial Lipschitz bound proves \eqref{regularized spacetime modulus}.

\end{proof}

We now combine the preceding local estimates with the
uniform boundary estimates of Corollary~\ref{bdrycont}
to remove the regularization parameters $\nu$ and $\zeta$.

\begin{prop}
\label{first compactness}
Under the assumptions of Theorem~\ref{existence}, let $a$ be as in
Lemma~\ref{bernstein}.
Fix $\epsilon\in(0,1)$ and $\delta\in[0,1]$.
For $\nu,\zeta\in(0,1)$, let $u_{\epsilon,\nu,\zeta}$ denote the
solutions to \eqref{pe} with $g_{\epsilon}$ replaced by
$\delta g_{\epsilon}$, constructed as in the proof of
Theorem~\ref{existence}.
Then there exists a unique viscosity solution
$u_{\epsilon}\in C(\overline{\Omega_{T}})$ to
\begin{align}
\label{pe1-delta}
\left\{\begin{aligned}
|Du_{\epsilon}|^{\kappa}F(D^{2}u_{\epsilon})
 -(u_{\epsilon})_{t}&=\delta g_{\epsilon}(u_{\epsilon})
 &&\text{in }\Omega_{T},\\
u_{\epsilon}&=\varphi+\epsilon
 &&\text{on }\partial_{p}\Omega_{T}.
\end{aligned}\right.
\end{align}
Moreover,
\begin{align}
\label{limit bound}
\epsilon\leq u_{\epsilon}\leq\|\varphi\|_{L^{\infty}(\partial_{p}\Omega_{T})}+\epsilon
\quad\text{in }\overline{\Omega_{T}},
\end{align}
and
\begin{align*}
u_{\epsilon,\nu,\zeta}\longrightarrow u_{\epsilon}
\quad\text{locally uniformly in }\Omega_{T}\text{ as }(\nu,\zeta)\to(0,0).
\end{align*}
If \eqref{uniform gradient quotient} holds for
$u=u_{\epsilon,\nu,\zeta}$ in a cylinder $Q\Subset\Omega_{T}$,
with the same $K$ for all $\nu,\zeta\in(0,1)$, then
\begin{align}
\label{power spatial lipschitz}
\left|u_{\epsilon}(x,t)^{1-\frac{a}{\kappa+2}}
 -u_{\epsilon}(y,t)^{1-\frac{a}{\kappa+2}}\right|
\leq\left(1-\frac{a}{\kappa+2}\right)
 K^{\frac{1}{\kappa+2}}|x-y|
\end{align}
for $[x,y]\times\{t\}\subset Q$.
\end{prop}
\begin{proof}
For every $B_{r}(x_{0})\Subset\Omega$ and $-T<s<t\leq0$, by Theorem~\ref{existence} and Lemma~\ref{bernstein}, we have \eqref{uniform gradient quotient} for $u=u_{\epsilon,\nu,\zeta}$ in $B_{r/2}(x_{0})\times[(s+t)/2,t]$, with $M=\|\varphi\|_{L^{\infty}(\partial_{p}\Omega_{T})}+\epsilon$ and $K$ independent of $\epsilon,\nu,\zeta$, and $\delta$.
By Theorem~\ref{existence} and Lemma~\ref{regularized time compactness},
the family $\{u_{\epsilon,\nu,\zeta}\}_{\nu,\zeta\in(0,1)}$ is
uniformly bounded and locally equicontinuous in $\Omega_{T}$.
By the Arzel\`a--Ascoli theorem, there exist a sequence
$(\nu_{j},\zeta_{j})\to(0,0)$ and a function $u_{\epsilon}$ such that
\begin{align*}
u_{\epsilon,\nu_{j},\zeta_{j}}\longrightarrow u_{\epsilon}
\quad\text{locally uniformly in }\Omega_{T}.
\end{align*}
By stability of viscosity solutions, $u_{\epsilon}$ satisfies the equation in \eqref{pe1-delta}.
Corollary~\ref{bdrycont} ensures that $u_{\epsilon}$ extends continuously to the parabolic boundary and satisfies $u_{\epsilon}=\varphi+\epsilon$ on $\partial_{p}\Omega_{T}$.
Passing to the limit in \eqref{es} yields \eqref{limit bound}.
Uniqueness follows from Lemma~\ref{comparison} and continuity at $t=0$.
Since every sequence with $(\nu_{j},\zeta_{j})\to(0,0)$ admits a locally uniformly convergent subsequence, the whole family converges locally uniformly to $u_{\epsilon}$ in $\Omega_{T}$.
Finally, by \eqref{uniform gradient quotient}, we have
\begin{align*}
\left|D\left(u_{\epsilon,\nu,\zeta}^{1-\frac{a}{\kappa+2}}\right)\right|
\leq\left(1-\frac{a}{\kappa+2}\right)
K^{\frac{1}{\kappa+2}}
\quad\text{in }Q.
\end{align*}
Integrating along $[x,y]\times\{t\}\subset Q$ and passing to the
limit proves \eqref{power spatial lipschitz}.
\end{proof}

\subsection{Uniform local estimates and weak spatial growth}
We next establish local estimates for bounded positive
viscosity solutions to \eqref{local cutoff equation}, uniformly in $\epsilon>0$
and $\delta\in[0,1]$.
After obtaining the spatial estimate, we derive a H\"older
estimate in time by applying a barrier argument to the
power $v_\epsilon^{2-a}$.

\begin{thm}
\label{equi conti}
Let $\kappa>0$, $0<\gamma<1$, and let $a$ be as in
Lemma~\ref{bernstein}. Let $F$ satisfy \ref{a1} and \ref{a2},
and let $v_{\epsilon}\in C(Q_{1})$ be a viscosity solution to
\begin{align}
\label{local cutoff equation}
|Dv_{\epsilon}|^{\kappa}F(D^{2}v_{\epsilon})
-(v_{\epsilon})_{t}=\delta g_{\epsilon}(v_{\epsilon})
\quad\text{in }Q_{1},
\end{align}
where $\epsilon>0$, $\delta\in[0,1]$, and
$0<v_{\epsilon}\leq M$ in $Q_{1}$.
Then, for every cylinder $Q'\Subset Q_{1}$,
\begin{align*}
\left|v_{\epsilon}(x,t)^{1-\frac{a}{\kappa+2}}
      -v_{\epsilon}(y,t)^{1-\frac{a}{\kappa+2}}\right|
\leq C|x-y|
\qquad\text{for }(x,t),(y,t)\in Q',\\
|v_{\epsilon}(x,t)-v_{\epsilon}(y,s)|
\leq C\bigl(|x-y|+|t-s|^{\frac{1}{2(2-a)}}\bigr)
\qquad\text{for }(x,t),(y,s)\in Q',
\end{align*}
where $C$ depends only on
$n,\lambda,\Lambda,\kappa,\gamma,M$, and $Q'$.
\end{thm}
\begin{proof}
Write $u=v_{\epsilon}$, fix a cylinder $Q'\Subset Q_{1}$, and choose
$Q=Q_{r}$ with $Q'\Subset Q\Subset Q_{1}$, depending only on $Q'$.
We first establish \eqref{power spatial lipschitz} for $u$ whenever $[x,y]\times\{t\}\subset Q'$.
For $\nu,\zeta\in(0,1)$, consider the regularized problem
\begin{align}
\label{local regularized problem}
\left\{\begin{aligned}
(\nu^{2}+|Dw_{\nu,\zeta}|^{2})^{\kappa/2}
F_{\zeta}(D^{2}w_{\nu,\zeta})-(w_{\nu,\zeta})_{t}
&=\delta g_{\epsilon}(w_{\nu,\zeta})&&\text{in }Q,\\
w_{\nu,\zeta}&=u&&\text{on }\partial_{p}Q.
\end{aligned}\right.
\end{align}
Let $c:=\min\bigl\{\epsilon,\min_{\partial_{p}Q}u\bigr\}>0$.
Since $g_{\epsilon}(c)=0$ and $g_{\epsilon}\geq0$, the constants $c$ and $M$ are respectively a subsolution and a supersolution to \eqref{local regularized problem}, with $c\leq u\leq M$ on $\partial_{p}Q$. 
Repeating the arguments in the proofs of Lemmas~\ref{bdrysuper}, \ref{bdrysub}, and Theorem~\ref{existence}, with $g_{\epsilon}$ replaced by $\delta g_{\epsilon}$, we obtain, for each $\nu,\zeta\in(0,1)$, a unique viscosity solution $w_{\nu,\zeta}\in C(\overline Q)\cap C^{\infty}(Q)$ to \eqref{local regularized problem} satisfying $c\leq w_{\nu,\zeta}\leq M$.
Lemma~\ref{bernstein} yields
\begin{align*}
\sup_{Q'}\frac{|Dw_{\nu,\zeta}|^{\kappa+2}}
{w_{\nu,\zeta}^{a}}\leq K,
\end{align*}
where $K>0$ depends only on $n,\lambda,\Lambda,\kappa,\gamma,M$, and $Q'$.
The arguments in the proofs of Lemma~\ref{regularized time compactness} and Proposition~\ref{first compactness} show that the family $\{w_{\nu,\zeta}\}$ is locally equicontinuous in $Q$ and that every locally uniform subsequential limit is a viscosity solution to \eqref{local cutoff equation} in $Q$ with boundary data $u$ on $\partial_{p}Q$.
Lemma~\ref{comparison} implies that every such limit coincides with $u$.
Hence
\begin{align*}
w_{\nu,\zeta}\longrightarrow u
\quad\text{locally uniformly in }Q
\text{ as }(\nu,\zeta)\to(0,0).
\end{align*}
The uniform gradient bound for
$w_{\nu,\zeta}^{1-a/(\kappa+2)}$ implies the corresponding
spatial Lipschitz estimate. Passing to the limit yields
\eqref{power spatial lipschitz} for $u$ in $Q'$.
Since the preceding argument applies to any interior cylinder, \eqref{power spatial lipschitz} also holds in $Q$, possibly with a larger $K$, and the first estimate follows.

Let $v:=u^{2-a}$.
Since $0<a<\gamma<1$, we have $2-a>1-a/(\kappa+2)>0$.
By the mean value theorem and \eqref{power spatial lipschitz},
using $0<u\leq M$, we obtain
\begin{align}
|v(x,t)-v(y,t)|
&=|u(x,t)^{2-a}-u(y,t)^{2-a}| \nonumber \\
&\leq\frac{2-a}{1-\frac{a}{\kappa+2}}
M^{1-a+\frac{a}{\kappa+2}}
\left|u(x,t)^{1-\frac{a}{\kappa+2}}
-u(y,t)^{1-\frac{a}{\kappa+2}}\right| \nonumber \\
&\leq\underbrace{(2-a)K^{\frac{1}{\kappa+2}}
M^{1-a+\frac{a}{\kappa+2}}}_{=:L}|x-y| \label{v lip bound}
\end{align}
for $(x,t),(y,t)\in Q$.
Let $\phi\in C^{2,1}$ touch $v$ from above or below at $z_{0}\in Q$.
Since $\phi(z_{0})=v(z_{0})=u(z_{0})^{2-a}>0$,  $\psi:=\phi^{1/(2-a)}$ is well-defined near $z_{0}$ and touches $u$ from the same side at $z_{0}$.
Then $\psi^{1-a/(\kappa+2)}$ touches $u^{1-a/(\kappa+2)}$ from the same side at $z_{0}$.
By \eqref{power spatial lipschitz} and Lemma~\ref{lipschitz test gradient}, we have
\begin{align*}
\left|D\left(\psi^{1-\frac{a}{\kappa+2}}\right)(z_{0})\right|
=\left(1-\frac{a}{\kappa+2}\right)u(z_{0})^{-\frac{a}{\kappa+2}}|D\psi(z_{0})|
\leq\left(1-\frac{a}{\kappa+2}\right)K^{\frac{1}{\kappa+2}},
\end{align*}
which implies
\begin{align}
\label{contact gradient bound}
u(z_{0})^{-a}|D\psi(z_{0})|^{\kappa+2}\leq K,
\qquad
|D\psi(z_{0})|^{\kappa}\leq(KM^{a})^{\frac{\kappa}{\kappa+2}}.
\end{align}
All quantities in the following computations are evaluated at $z_{0}$.
By direct calculation, we have
\begin{align}
\label{phi diff}
\phi_{t}&=(2-a)u^{1-a}\psi_{t}, \qquad D^{2}\phi=(2-a)u^{1-a}D^{2}\psi
 +(2-a)(1-a)u^{-a}D\psi\otimes D\psi.
\end{align}
Note that by \ref{a1} and $0<a<1$, we have
\begin{align}
\M^{+}(D^{2}\phi)
&\geq \M^{+}((2-a)u^{1-a}D^{2}\psi)+\M^{-}((2-a)(1-a)u^{-a}D\psi\otimes D\psi) \nonumber \\
&\geq (2-a)u^{1-a}F(D^{2}\psi) +\lambda(2-a)(1-a)u^{-a}|D\psi|^{2} \label{pucci lower}
\end{align}
and
\begin{align}
\M^{-}(D^{2}\phi)
&\leq \M^{-}((2-a)u^{1-a}D^{2}\psi)+\M^{+}((2-a)(1-a)u^{-a}D\psi\otimes D\psi) \nonumber \\
&\leq (2-a)u^{1-a}F(D^{2}\psi) +\Lambda(2-a)(1-a)u^{-a}|D\psi|^{2}.  \label{pucci upper} 
\end{align}
If $\phi$ touches $v$ from above at $z_{0}$, then $\psi$ touches $u$ from above at $z_{0}$, and by \eqref{local cutoff equation} this yields
\begin{align}
\label{sub ineq}
\psi_{t}-|D\psi|^{\kappa}F(D^{2}\psi)\leq -\delta\beta_{\epsilon}(u)u^{\gamma-1}.
\end{align}
By \eqref{phi diff}, \eqref{pucci lower}, and \eqref{sub ineq}, we obtain
\begin{align}
\phi_{t}-|D\psi|^{\kappa}\M^{+}(D^{2}\phi)
&\leq (2-a)u^{1-a}\psi_t-(2-a)u^{1-a}|D\psi|^{\kappa}F(D^{2}\psi) -\lambda(2-a)(1-a)u^{-a}|D\psi|^{\kappa+2} \nonumber \\
&= (2-a)u^{1-a}
\left(\psi_t-|D\psi|^{\kappa}F(D^{2}\psi)\right)
-\lambda(2-a)(1-a)u^{-a}|D\psi|^{\kappa+2} \nonumber \\
&\leq -(2-a)\delta\beta_{\epsilon}(u)u^{\gamma-a}-\lambda(2-a)(1-a)u^{-a}|D\psi|^{\kappa+2} \label{time diff upper}
\leq0.
\end{align}
If $\phi$ touches $v$ from below at $z_{0}$, then $\psi$ touches $u$ from below at $z_{0}$, and by \eqref{local cutoff equation} this yields
\begin{align}
\label{sup ineq}
\psi_{t}-|D\psi|^{\kappa}F(D^{2}\psi)\geq -\delta\beta_{\epsilon}(u)u^{\gamma-1}.
\end{align}
Similarly, by \eqref{phi diff}, \eqref{pucci upper}, \eqref{sup ineq}, and \eqref{contact gradient bound}, we get
\begin{align}
\phi_{t}-|D\psi|^{\kappa}\M^{-}(D^{2}\phi)
&\geq (2-a)u^{1-a}\psi_t-(2-a)u^{1-a}|D\psi|^{\kappa}F(D^{2}\psi) -\Lambda(2-a)(1-a)u^{-a}|D\psi|^{\kappa+2} \nonumber \\
&= (2-a)u^{1-a}
\left(\psi_t-|D\psi|^{\kappa}F(D^{2}\psi)\right)
-\Lambda(2-a)(1-a)u^{-a}|D\psi|^{\kappa+2} \nonumber \\
&\geq -(2-a)\delta\beta_{\epsilon}(u)u^{\gamma-a}-\Lambda(2-a)(1-a)u^{-a}|D\psi|^{\kappa+2} \nonumber \\
&\geq -(2-a)M^{\gamma-a}-\Lambda(2-a)(1-a)K. \label{time diff lower}
\end{align}
Choose $C_{1}\geq1$ such that
\begin{align*}
C_{1}\geq2n\Lambda(KM^{a})^{\frac{\kappa}{\kappa+2}}, \qquad
C_{1}\geq\Lambda(2-a)(1-a)K+(2-a)M^{\gamma-a}.
\end{align*}
Hence \eqref{time diff upper}, \eqref{time diff lower}, and
\eqref{contact gradient bound} imply, for $A\geq0$,
\begin{align}
\label{quadratic contacts}
\begin{aligned}
D^{2}\phi=2AI
 &\quad\Longrightarrow\quad\phi_{t}\leq C_{1}A
 &&\text{if $\phi$ touches $v$ from above},\\
D^{2}\phi=-2AI
 &\quad\Longrightarrow\quad\phi_{t}\geq-C_{1}(A+1)
 &&\text{if $\phi$ touches $v$ from below}.
\end{aligned}
\end{align}
Fix $(x,s),(x,t)\in Q'$ with $s<t$, and choose $R>0$,
depending only on $Q',Q$, such that
$\overline{B_{R}(x)}\times[s,t]\subset Q$.
For $\eta>0$, define
\begin{align*}
A_{\eta}:=\frac{L^{2}}{4\eta}+\frac{M^{2-a}}{R^{2}}, \qquad
h^{\pm}(y,\tau):=v(x,s)\pm\left\{
\eta+A_{\eta}|y-x|^{2}+2C_{1}(A_{\eta}+1)(\tau-s)\right\}.
\end{align*}
We claim that $h^{-}\leq v \leq h^{+}$ in $\overline{B_{R}(x)}\times[s,t]$.
As in the proof of Lemma~\ref{regularized time compactness}, \eqref{v lip bound}, $0<v\leq M^{2-a}$, and $A_{\eta}R^{2}\geq M^{2-a}$ imply $h^{-}\leq v\leq h^{+}$ on $\partial_{p}(B_{R}(x)\times(s,t])$.
We now assume that $v-h^{+}$ attains a positive maximum at $z_{0}=(y_{0},\tau_{0})$
in $\overline{B_{R}(x)}\times[s,t]$.
The boundary inequalities imply that $z_{0}\in B_{R}(x)\times(s,t]$.
Then
\begin{align*}
\phi(y,\tau):=h^{+}(y,\tau)+v(z_{0})-h^{+}(z_{0})
\end{align*}
touches $v$ from above at $z_{0}$, with $D^{2}\phi=2A_{\eta}I$ and $\phi_{t}=2C_{1}(A_{\eta}+1)$.
By \eqref{quadratic contacts}, we have at $z_{0}$
\begin{align*}
2C_{1}(A_{\eta}+1)=\phi_{t}\leq C_{1}A_{\eta},
\end{align*}
which is a contradiction.
Similarly, assume that $h^{-}-v$ attains a positive maximum at $z_{1}$
in $\overline{B_{R}(x)}\times[s,t]$.
Then $z_{1}\in B_{R}(x)\times(s,t]$ and
$h^{-}+v(z_{1})-h^{-}(z_{1})$ touches $v$ from below at $z_{1}$.
By \eqref{quadratic contacts}, we obtain
\begin{align*}
-2C_{1}(A_{\eta}+1)=h^{-}_{t}\geq-C_{1}(A_{\eta}+1),
\end{align*}
which is again a contradiction.
Hence $h^{-}\leq v\leq h^{+}$ in
$\overline{B_{R}(x)}\times[s,t]$, and
\begin{align*}
|v(x,t)-v(x,s)|
\leq\eta+2C_{1}\left(\frac{L^{2}}{4\eta}
 +\frac{M^{2-a}}{R^{2}}+1\right)(t-s).
\end{align*}
Since $0<t-s\leq1$, taking $\eta=\sqrt{t-s}$ yields
\begin{align}
\label{power time modulus}
|v(x,t)-v(x,s)|\leq C|t-s|^{1/2}.
\end{align}
Using $|A^{q}-B^{q}|\leq|A-B|^{q}$ for $A,B\geq0$ and $0<q\leq1$, with $q=1/(2-a)$, we have
\begin{align*}
|u(x,t)-u(x,s)|
\leq|v(x,t)-v(x,s)|^{\frac{1}{2-a}}
\leq C|t-s|^{\frac{1}{2(2-a)}}.
\end{align*}
By the mean value theorem and \eqref{power spatial lipschitz}, we also have
\begin{align*}
|u(x,t)-u(y,t)|\leq \frac{1}{1-\frac{a}{\kappa+2}}
M^{\frac{a}{\kappa+2}}
\left|u(x,t)^{1-\frac{a}{\kappa+2}}-u(y,t)^{1-\frac{a}{\kappa+2}}\right|    
\leq(KM^{a})^{\frac{1}{\kappa+2}}|x-y|
\quad\text{for }(x,t),(y,t)\in Q'.
\end{align*}
Combining the above estimates proves the second estimate.
\end{proof}

The preceding estimates pass to locally uniform limits
and yield superlinear spatial growth at zeros
of the limiting functions.
This will be used in Section~\ref{sec5} to identify the
homogeneous equation across the zero set.

\begin{cor}
\label{weak growth corollary}
Let $\kappa>0$ and $0<\gamma<1$, and let $a$ be as in
Lemma~\ref{bernstein}.
Let $\mathcal V$ be a locally uniformly bounded family of positive
viscosity solutions as in Theorem~\ref{equi conti}, with common ellipticity
constants $0<\lambda\leq\Lambda$.
Then $\mathcal V$ is locally equicontinuous in $Q_{1}$.
If $V$ is a locally uniform limit of a sequence in $\mathcal V$, then,
for every $0<r<R<1$,
\begin{align}
\label{weak growth}
V(x_{0},t_{0})=0
\quad\Longrightarrow\quad
0\leq V(x,t_{0})\leq C|x-x_{0}|^{\frac{\kappa+2}{\kappa+2-a}}
\end{align}
for all $(x_{0},t_{0})\in Q_{r}$ and $x\in B_{r}$, where
$C=C\Bigl(n,\lambda,\Lambda,\kappa,\gamma,r,R, \sup_{v\in\mathcal V}\|v\|_{L^{\infty}(Q_R)}\Bigr)>0$.
In particular, $DV=0$ on $\{V=0\}\cap Q_1$.
\end{cor}
\begin{proof}
By Theorem~\ref{equi conti} and scaling, $\mathcal V$ is locally equicontinuous in $Q_1$ and satisfies the corresponding spatial estimate uniformly on interior cylinders.
Passing to the locally uniform limit, for $(x_0,t_0)\in Q_r$ with $V(x_0,t_0)=0$, we obtain
\begin{align*}
V(x,t_0)^{1-\frac{a}{\kappa+2}}\leq C|x-x_0|
\qquad\text{for }x\in B_r.
\end{align*}
This proves \eqref{weak growth}. Since
$(\kappa+2)/(\kappa+2-a)>1$, we conclude that $DV=0$ on $\{V=0\}\cap Q_1$.
\end{proof}

\subsection{Existence for the free boundary problem}
We conclude this section by letting $\epsilon\to0$ in \eqref{pe1}.
Comparison gives monotonicity in $\epsilon$, while the
local estimates and boundary barriers ensure local uniform
convergence and attainment of the prescribed boundary data.

\begin{prop}
\label{existence and equiconti}
Under the assumptions of Theorem~\ref{existence}, there exists
a viscosity solution $u\in C(\ov{\Omega_{T}})$ to \eqref{eq}
satisfying
\begin{align}
\label{u bound}
0\leq u\leq\|\varphi\|_{L^{\infty}(\partial_{p}\Omega_{T})}
\quad\text{in }\overline{\Omega_{T}}.
\end{align}
Moreover, the solutions $u_{\epsilon}$ to \eqref{pe1} given by
Proposition~\ref{first compactness} converge locally uniformly to $u$
in $\Omega_{T}$ as $\epsilon\to0^{+}$.
\end{prop}
\begin{proof}
Note that for $0<\epsilon_{1}<\epsilon_{2}<1$, by the definition of $g_{\epsilon}$, we have $g_{\epsilon_{1}}\geq g_{\epsilon_{2}}$, which implies that
$u_{\epsilon_{1}}$ satisfies 
\begin{align*}
|Du_{\epsilon_{1}}|^{\kappa}F(D^{2}u_{\epsilon_{1}})-(u_{\epsilon_{1}})_{t}\geq g_{\epsilon_{2}}(u_{\epsilon_{1}})
\quad\text{in }\Omega_{T},
\end{align*}
in the viscosity sense.
Since $u_{\epsilon_{1}}=\varphi+\epsilon_{1}< \varphi+\epsilon_{2}=u_{\epsilon_{2}}$ on $\partial_{p}\Omega_{T}$, by Lemma~\ref{comparison}, continuity at $t=0$, and \eqref{limit bound},
we obtain
\begin{align*}
0< u_{\epsilon_{1}}\leq u_{\epsilon_{2}}
 \leq\|\varphi\|_{L^{\infty}(\partial_{p}\Omega_{T})}+\epsilon_{2}< \norm{\varphi}_{L^{\infty}(\partial_{p}\Omega_{T})}+1
\quad\text{in }\overline{\Omega_{T}}.
\end{align*}
Consequently, the pointwise limit $u:=\lim_{\epsilon\to0^{+}}u_{\epsilon}$ exists in $\Omega_{T}$.
By Theorem~\ref{equi conti} and the Arzel\`a--Ascoli theorem,
every sequence $\{u_{\epsilon_j}\}$ with $\epsilon_j\to0^{+}$ admits a locally uniformly convergent subsequence, whose limit coincides with $u$.
Thus $u\in C(\Omega_{T})$ and $u_{\epsilon}\to u$ locally uniformly in $\Omega_{T}$.
Moreover, Corollary~\ref{bdrycont} and Proposition~\ref{first compactness} imply that $u\in C(\ov{\Omega_{T}})$ with $u=\varphi$ on $\partial_{p}\Omega_{T}$.
Passing to the limit in \eqref{limit bound} yields \eqref{u bound}.
It remains to prove that $u$ satisfies the equation in \eqref{eq}. Let $Q\Subset\{u>0\}\cap\Omega_{T}$.
From the uniform convergence of $u_{\epsilon}$ to $u$ on $\ov{Q}$, for all sufficiently small $\epsilon$, we have
\begin{align*}
u_{\epsilon}\geq (\inf_{\ov{Q}}u)/2>2\epsilon \quad \text{in } Q,\qquad
\beta_{\epsilon}(u_{\epsilon})=1
\quad\text{in }Q,
\end{align*}
and hence
\begin{align*}
|Du_{\epsilon}|^{\kappa}F(D^{2}u_{\epsilon})-(u_{\epsilon})_{t}=u_{\epsilon}^{\gamma-1} \quad \text{in } Q.   
\end{align*}
By stability of viscosity solutions, $u$ satisfies the equation
in \eqref{eq}.
\end{proof}

\section{Improvement and optimal growth along the free boundary}
\label{sec5}
We first prove the optimal upper growth estimate in backward cylinders by compactness and iteration.
The main step is to identify the homogeneous equation across the zero set of a limit as the absorption coefficient tends to zero.
We consider positive viscosity solutions to
\begin{align}
\label{pe2 eq}
|Dv|^{\kappa}F(D^{2}v)-v_{t}=\delta g_{\epsilon}(v),
\end{align}
where $\epsilon>0$ and $\delta\geq0$.

\subsection{A comparison function and the homogeneous limit}
We combine the spatial profile from Lemma \ref{barrier homo}
with a decreasing positive function of time to construct
a subsolution to \eqref{pe2 eq}.

\begin{lem}
\label{comparison bump}
Let $\kappa>0$, $0<\gamma<1$, and let $F$ satisfy \ref{a1}.
Let $\phi$, $K_{0}$, and $k_{0}$ be as in Lemma~\ref{barrier homo}.
Let $x_{0}\in\R^{n}$, $s<t_{1}$, $A_{0}>0$, and $r>0$. Define
\begin{align*}
A(t)&:=\left[A_{0}^{-\kappa}+\kappa K_{0}r^{-(\kappa+2)}(t-s)\right]^{-1/\kappa},
\qquad t\in[s,t_{1}],\\
W(x,t)&:=A(t)\phi\left(\frac{x-x_{0}}{r}\right),
\qquad (x,t)\in\ov{B_{r}(x_{0})}\times[s,t_{1}].
\end{align*}
If
\begin{align}
\label{delta range}
0\leq\delta\leq k_{0}r^{-(\kappa+2)}A(t_{1})^{\kappa+2-\gamma},
\end{align}
then, for every $\epsilon>0$, $W$ is a viscosity subsolution to
\eqref{pe2 eq} in $B_{r}(x_{0})\times(s,t_{1})$, with $W=0$ on
$\partial B_{r}(x_{0})\times[s,t_{1}]$.
\end{lem}

\begin{proof}
Note that $W$ is smooth and positive in
$B_{r}(x_{0})\times(s,t_{1})$.
By direct calculation, we have
\begin{align*}
DW(x,t)&=\frac{A(t)}{r}D\phi\left(\frac{x-x_{0}}{r}\right),\\
D^{2}W(x,t)&=\frac{A(t)}{r^{2}}D^{2}\phi\left(\frac{x-x_{0}}{r}\right),\\
W_{t}(x,t)&=A'(t)\phi\left(\frac{x-x_{0}}{r}\right).
\end{align*}
Moreover, $A$ is decreasing and satisfies 
\begin{align}
\label{A rel}
A'=-K_{0}r^{-(\kappa+2)}A^{\kappa+1}.    
\end{align}
By \ref{a1}, \eqref{A rel}, Lemma \ref{barrier homo} with $\nu=0$, and \eqref{delta range}, we obtain
\begin{align*}
|DW|^{\kappa}F(D^{2}W)-W_{t}&\geq |DW|^{\kappa}\M^{-}(D^{2}W)-W_{t}\\
&= r^{-(\kappa+2)}A^{\kappa+1}
\left(|D\phi|^{\kappa}\M^{-}(D^{2}\phi)+K_{0}\phi\right)\\
&\geq k_{0}r^{-(\kappa+2)}A^{\kappa+2-\gamma}W^{\gamma-1}\\
&\geq k_{0}r^{-(\kappa+2)}A(t_{1})^{\kappa+2-\gamma}W^{\gamma-1}\\
&\geq\delta g_{\epsilon}(W) \qquad \text{in } B_{r}(x_{0})\times(s,t_{1}).
\end{align*}
Finally, $W=0$ on $\partial B_{r}(x_{0})\times[s,t_{1}]$ since $\phi=0$ on $\partial B_{1}$.
\end{proof}

Comparison with this subsolution gives persistence of
positivity for locally uniform limits as the absorption
coefficient tends to zero.

\begin{lem}
\label{positivity}
Let $\kappa>0$, $0<\gamma<1$, and
$D=U\times(\tau_{0},\tau_{1}]$, where $U\subset\R^{n}$ is open and $\tau_{0}<\tau_{1}$.
For each $j$, let $F_{j}$ satisfy \ref{a1} with the same ellipticity
constants $0<\lambda\leq\Lambda$, and let $v_{j}\in C(D)$ be a
positive viscosity solution to
\begin{align*}
|Dv_{j}|^{\kappa}F_{j}(D^{2}v_{j})-(v_{j})_{t}
=\delta_{j}g_{\epsilon_{j}}(v_{j})\qquad\text{in }D,
\end{align*}
where $\epsilon_{j}>0$ and $\delta_{j}\geq0$.
Assume that $\delta_{j}\to0$ and that $v_{j}\to V$ locally uniformly
in $D$ as $j\to\infty$.
If $s<t_{1}$, $B_{r_{0}}(x_{0})\times[s,t_{1}]\Subset D$, and
$V(x_{0},s)>0$, then $V(x_{0},t_{1})>0$.
\end{lem}

\begin{proof}
By continuity, there exists $r\in(0,r_{0})$ such that
\begin{align*}
V(x,s)\geq\frac{V(x_{0},s)}{2}
\qquad\text{for }x\in\overline{B_{r}(x_{0})}.
\end{align*}
Let $W$ be the function in Lemma~\ref{comparison bump} with
$A_{0}=V(x_{0},s)/4$. 
Note that $W$ is independent of $j$.
Since $A(t_{1})>0$ and $\delta_{j}\to0$, for all sufficiently large
$j$ we have
\begin{align*}
\delta_{j}\leq k_{0}r^{-(\kappa+2)}A(t_{1})^{\kappa+2-\gamma}.
\end{align*}
By local uniform convergence and the positivity of $v_{j}$, we also
have, for all sufficiently large $j$,
\begin{align*}
W(x,s)&\leq A_{0}\leq v_{j}(x,s)
\qquad\text{on }\overline{B_{r}(x_{0})},\\
W&=0\leq v_{j}
\qquad\text{on }\partial B_{r}(x_{0})\times[s,t_{1}].
\end{align*}
For each such $j$, Lemma~\ref{comparison bump}, applied with
$F=F_{j}$, $\delta=\delta_{j}$, and $\epsilon=\epsilon_{j}$, shows
that $W$ is a viscosity subsolution to the equation satisfied by $v_{j}$. 
By Lemma~\ref{comparison} and continuity, we have
\begin{align*}
W\leq v_{j}\qquad\text{in }B_{r}(x_{0})\times[s,t_{1}].
\end{align*}
Letting $j\to\infty$, we obtain
\begin{align*}
V(x_{0},t_{1})\geq W(x_{0},t_{1})=A(t_{1})>0.
\end{align*}
\end{proof}

We combine this persistence property with the weak spatial
growth estimate of Corollary~\ref{weak growth corollary}
to identify the homogeneous equation across the zero set
of the limit.

\begin{lem}
\label{homogeneous limit}
Let $\kappa>0$, $0<\gamma<1$, and
$D=U\times(\tau_{0},\tau_{1}]$, where $U\subset\R^{n}$ is open
and $\tau_{0}<\tau_{1}$.
For each $j$, let $F_{j}$ satisfy \ref{a1} and \ref{a2}
with the same ellipticity constants $0<\lambda\leq\Lambda$,
and let $v_{j}\in C(D)$ be a positive viscosity solution to
\begin{align*}
|Dv_{j}|^{\kappa}F_{j}(D^{2}v_{j})-(v_{j})_{t}
=\delta_{j}g_{\epsilon_{j}}(v_{j})
\qquad\text{in }D,
\end{align*}
where $\epsilon_{j}>0$ and $\delta_{j}\geq0$.
Assume that $\delta_{j}\to0$ and that $0<v_{j}\leq M$ in $D$
for some $M>0$ independent of $j$.
Then, after passing to a subsequence, $v_{j}$ converges locally
uniformly in $D$ to a function $V\in C(D)$ with $0\leq V\leq M$.
Moreover, $V$ is a viscosity solution to
\begin{align}
\label{homogeneous equation}
|DV|^{\kappa}F_{\infty}(D^{2}V)-V_{t}=0
\qquad\text{in }D,
\end{align}
for some $F_{\infty}$ satisfying \ref{a1} and \ref{a2}.
\end{lem}

\begin{proof}
Since $0\leq\delta_{j}\leq1$ for all sufficiently large $j$,
Theorem~\ref{equi conti} implies that,
after passing to a subsequence,
\begin{align*}
v_{j}\to V\quad\text{locally uniformly in }D,
\qquad V\in C(D),\qquad 0\leq V\leq M.
\end{align*}
By \ref{a1}, after passing to a further subsequence,
$F_{j}\to F_{\infty}$ locally uniformly on $\mathcal{S}^{n}$,
where $F_{\infty}$ satisfies \ref{a1} and \ref{a2}.

Since $\delta_{j}g_{\epsilon_{j}}\geq0$, each $v_{j}$ satisfies
\begin{align*}
|Dv_{j}|^{\kappa}F_{j}(D^{2}v_{j})-(v_{j})_{t}\geq0
\qquad\text{in }D,
\end{align*}
in the viscosity sense. By stability of viscosity solutions, $V$ is a
viscosity subsolution to \eqref{homogeneous equation} in $D$.

Fix $z_{0}=(x_{0},t_{0})\in\{V>0\}\cap D$.
By continuity and local uniform convergence, there exists a
neighborhood $N$ of $z_{0}$ relative to $D$ such that, for all
sufficiently large $j$, $v_{j}\geq V(z_{0})/2>0$ in $N$.
Since 
\begin{align*}
0\leq\delta_{j}g_{\epsilon_{j}}(v_{j})
\leq\delta_{j}\left(V(z_{0})/2\right)^{\gamma-1}
\to0\qquad\text{uniformly in }N,
\end{align*}
stability of viscosity solutions now yields \eqref{homogeneous equation} in $N$.
Thus $V$ satisfies \eqref{homogeneous equation} in
$\{V>0\}\cap D$.

It remains to prove the supersolution inequality at a point
$z_{0}=(x_{0},t_{0})\in\{V=0\}\cap D$.
Let $\xi\in C^{2,1}$ touch $V$ from below at $z_{0}$, with
$\xi(z_{0})=V(z_{0})=0$.
By Corollary~\ref{weak growth corollary}, $V(\cdot,t_{0})$ is differentiable at $x_{0}$ and
$DV(z_{0})=0$. Since $\xi$ touches $V$ from below at $z_{0}$, we obtain
$D\xi(z_{0})=0$.

Choose any $r,h_{0}>0$ such that
\begin{align*}
B_{r}(x_{0})\times[t_{0}-h_{0},t_{0}]\Subset D.
\end{align*}
If $V(x_{0},t)>0$ for some $t\in[t_{0}-h_{0},t_{0})$,
Lemma~\ref{positivity}, applied between $t$ and $t_{0}$, would give
$V(x_{0},t_{0})>0$, a contradiction. Hence
\begin{align*}
V(x_{0},t)=0\qquad\text{for every }t\in[t_{0}-h_{0},t_{0}].
\end{align*}
For all sufficiently small $h\in(0,h_{0})$, Taylor expansion in time therefore yields
\begin{align*}
-h\xi_{t}(z_{0})+o(h)
=\xi(x_{0},t_{0}-h)-\xi(z_{0})
\leq V(x_{0},t_{0}-h)-V(z_{0})=0.
\end{align*}
Dividing by $h$ and letting $h\to0^{+}$, we obtain
$\xi_{t}(z_{0})\geq0$.
Since $\kappa>0$ and $D\xi(z_{0})=0$, it follows that
\begin{align*}
|D\xi(z_{0})|^{\kappa}F_{\infty}(D^{2}\xi(z_{0}))
-\xi_{t}(z_{0})=-\xi_{t}(z_{0})\leq0.
\end{align*}
This proves the supersolution inequality at $z_{0}$ and
completes the proof.
\end{proof}

\subsection{One-step improvement and iteration}
We combine Lemma~\ref{homogeneous limit} with the optimal
spatial regularity estimate in Lemma~\ref{time derivative}.
Since $(\kappa+2)/(\kappa+1)>\alpha$, a compactness argument
gives the following improvement at a fixed scale.

\begin{lem}
\label{one step}
Let $F$ satisfy \ref{a1} and \ref{a2}, and let $v$ be a
positive viscosity solution to \eqref{pe2 eq} in $Q_{1}^{\theta}$,
where $\epsilon>0$ and $\delta\geq0$.
There exist constants $\rho\in(0,1/4)$ and
$\delta_{0},\eta_{0}\in(0,1)$, depending only on
$n,\lambda,\Lambda,\kappa$, and $\gamma$, such that, if
\begin{align*}
0<v\leq1\quad\text{in }Q_{1}^{\theta},\qquad
0\leq\delta\leq\delta_{0},\qquad
v(0,0)\leq\eta_{0},
\end{align*}
then
\begin{align}
\label{one step estimate}
\sup_{Q_{\rho}^{\theta}}v\leq\rho^{\alpha}.
\end{align}
\end{lem}

\begin{proof}
Let $C_{1}=C_{1}(n,\lambda,\Lambda,\kappa)\geq1$ be the
constant in the spatial gradient estimate of
Lemma~\ref{time derivative} for solutions bounded by $1$.
Since $\alpha<(\kappa+2)/(\kappa+1)$, choose
$\rho\in(0,1/4)$ such that
\begin{align*}
C_{1}\rho^{(\kappa+2)/(\kappa+1)}<\rho^{\alpha}/2.
\end{align*}
Suppose, to the contrary, that the conclusion fails for this $\rho$.
Then, for every $j\in\mathbb N$, there exist $\epsilon_{j}>0$,
$\delta_{j}\geq0$, $F_{j}$ satisfying \ref{a1} and
\ref{a2} with the same ellipticity constants $0<\lambda\leq\Lambda$, and a positive
viscosity solution $v_{j}$ to
\begin{align*}
|Dv_{j}|^{\kappa}F_{j}(D^{2}v_{j})-(v_{j})_{t}
=\delta_{j}g_{\epsilon_{j}}(v_{j})
\quad\text{in }Q_{1}^{\theta},
\end{align*}
such that
\begin{align*}
0<v_{j}\leq1\quad\text{in }Q_{1}^{\theta},\qquad
0\leq\delta_{j}\leq1/j,\qquad
v_{j}(0,0)\leq1/j,
\end{align*}
but
\begin{align}
\label{contra}
\sup_{Q_{\rho}^{\theta}}v_{j}>\rho^{\alpha}.
\end{align}
By Lemma~\ref{homogeneous limit}, after passing to a subsequence, 
$v_{j}$ converges locally uniformly in $Q_{1}^{\theta}$ to
a function $V\in C(Q_{1}^{\theta})$ that satisfies
\eqref{homogeneous equation} in $Q_{1}^{\theta}$ for some
$F_{\infty}$ satisfying \ref{a1} and \ref{a2}.
Moreover,
\begin{align*}
0\leq V\leq1\quad\text{in }Q_{1}^{\theta},
\qquad V(0,0)=0.
\end{align*}
By Lemma~\ref{positivity}, we obtain
\begin{align*}
V(0,t)=0\qquad\text{for }-1<t\leq0.
\end{align*}
Lemma~\ref{time derivative} yields
\begin{align*}
|DV(x,t)-DV(y,t)|\leq C_{1}|x-y|^{1/(\kappa+1)}
\end{align*}
for $(x,t),(y,t)\in\ov{Q_{1/2}}$.
Nonnegativity of $V$ and $V(0,t)=0$ imply $DV(0,t)=0$ for
$-1/4\leq t\leq0$. Hence, for $(x,t)\in\ov{Q_{1/2}}$,
\begin{align}
V(x,t)=V(x,t)-V(0,t)
&=\int_{0}^{1}DV(\tau x,t)\cdot x\,d\tau\nonumber\\
&=\int_{0}^{1}\bigl(DV(\tau x,t)-DV(0,t)\bigr)\cdot x\,d\tau\nonumber\\
&\leq C_{1}|x|^{(\kappa+2)/(\kappa+1)}
\int_{0}^{1}\tau^{1/(\kappa+1)}\,d\tau\nonumber\\
&\leq C_{1}|x|^{(\kappa+2)/(\kappa+1)}.\label{C1}
\end{align}
Since $\theta>1$ and $\rho<1/4$, we have
$\ov{Q_{\rho}^{\theta}}\subset\ov{Q_{1/2}}$.
Our choice of $\rho$ therefore gives
\begin{align*}
\sup_{Q_{\rho}^{\theta}}V
\leq C_{1}\rho^{(\kappa+2)/(\kappa+1)}
<\rho^{\alpha}/2.
\end{align*}
By uniform convergence on $\ov{Q_{\rho}^{\theta}}$, we have
$\sup_{Q_{\rho}^{\theta}}v_{j}\leq\rho^{\alpha}$
for all sufficiently large $j$, contradicting \eqref{contra}.
This completes the proof.
\end{proof}

We iterate Lemma~\ref{one step} as long as the rescaled
value at the center satisfies the required smallness
condition.
At the remaining scales, the value at the center controls
the supremum.

\begin{lem}
\label{iterated growth}
Let $F$ satisfy \ref{a1} and \ref{a2}, and let $v$ be a
positive viscosity solution to
\eqref{pe2 eq} in $Q_{1}^{\theta}$,
where $\epsilon>0$ and $\delta\geq0$. Assume that $0<v\leq1$ in
$Q_{1}^{\theta}$.
There exist constants $C\geq1$ and $\delta_{0}\in(0,1)$, depending only on
$n,\lambda,\Lambda,\kappa$, and $\gamma$, such that, if
$0\leq\delta\leq\delta_{0}$, then
\begin{align*}
\sup_{Q_{r}^{\theta}}v\leq C\bigl(v(0,0)+r^{\alpha}\bigr)
\qquad\text{for }0<r\leq1.
\end{align*}
\end{lem}

\begin{proof}
Let $\rho\in(0, 1/4)$ and $\delta_{0}, \eta_{0}\in(0,1)$ be the constants in Lemma \ref{one step}.
If $v(0,0)>\eta_{0}$, then
\begin{align*}
\sup_{Q_{r}^{\theta}}v\leq1<\eta_{0}^{-1}v(0,0),
\end{align*}
so the assertion follows. Suppose $v(0,0)\leq\eta_{0}$, and let
$N\geq0$ be the largest integer such that $v(0,0)\leq\eta_{0}\rho^{N\alpha}$.
For $0\leq j\leq N$, define $v_{j}$ by
\begin{align*}
v_{j}(x,t):=\rho^{-j\alpha}v(\rho^{j}x,\rho^{j\theta}t) \quad \text{in }Q_{1}^{\theta}.
\end{align*}
Then $v_{j}$ satisfies
\begin{align*}
|Dv_{j}|^{\kappa}F_{j}(D^{2}v_{j})-(v_{j})_{t}
=\delta g_{\epsilon/\rho^{j\alpha}}(v_{j}) \quad \text{in } Q_{1}^{\theta},
\end{align*}
in the viscosity sense, where
\begin{align*}
F_{j}(X):=\rho^{j(2-\alpha)}F(\rho^{j(\alpha-2)}X)
\quad\text{for }X\in\mathcal{S}^{n}.
\end{align*}
Here we used $\kappa(1-\alpha)+2-\alpha=\theta-\alpha=\alpha(1-\gamma)$,
and we note that each $F_{j}$ satisfies \ref{a1} and \ref{a2} with the same ellipticity constants.
Moreover, we have
\begin{align*}
v_{j}(0,0)=\rho^{-j\alpha}v(0,0)
\leq\eta_{0}\rho^{(N-j)\alpha}\leq\eta_{0}.
\end{align*}
Starting with $v_{0}=v$, Lemma~\ref{one step} applies
inductively, since the preceding step ensures
$0<v_{j}\leq1$ in $Q_{1}^{\theta}$.
Consequently,
\begin{align*}
\sup_{Q_{\rho^{j}}^{\theta}}v\leq\rho^{j\alpha}
\qquad\text{for }0\leq j\leq N+1.
\end{align*}
If $\rho^{j+1}<r\leq\rho^{j}$ for some $0\leq j\leq N$, then
\begin{align*}
\sup_{Q_{r}^{\theta}}v\leq\rho^{j\alpha}
\leq\rho^{-\alpha}r^{\alpha}.
\end{align*}
If $0<r\leq\rho^{N+1}$, the maximality of $N$ gives
\begin{align*}
\sup_{Q_{r}^{\theta}}v\leq\rho^{(N+1)\alpha}
<\eta_{0}^{-1}v(0,0).
\end{align*}
Choosing $C=\max\left\{\rho^{-\alpha}, \eta_{0}^{-1}\right\}$ completes the proof.
\end{proof}

By rescaling and applying Lemma~\ref{iterated growth}, we obtain the following estimate for $\delta=1$.

\begin{thm}
\label{optimal growth proposition}
Let $F$ satisfy \ref{a1} and \ref{a2}. 
Let $v$ be a positive viscosity solution to \eqref{pe2 eq} in $Q_{R}^{\theta}(z_{0})$,
with $\epsilon>0$ and $\delta=1$. Assume that $0<v\leq M$ in
$Q_{R}^{\theta}(z_{0})$. Then
\begin{align}
\label{optimal growth}
\sup_{Q_{r}^{\theta}(z_{0})}v
\leq C\left\{v(z_{0})+
\left[1+\left(\frac{M}{R^{\alpha}}\right)^{\frac{\kappa+2-\gamma}{2-\gamma}}\right]r^{\alpha}\right\}
\qquad \text{for } 0<r\leq R,
\end{align}
where $C=C(n,\lambda,\Lambda,\kappa,\gamma)$.
\end{thm}

\begin{proof}
Let $\delta_{0}$ be the constant in Lemma \ref{iterated growth}, and write $z_{0}=(x_{0},t_{0})$.
Define
\begin{align*}
K:=\max\left\{1,\delta_{0}^{-1/(\kappa+2-\gamma)},
\frac{M}{R^{\alpha}}\right\},\qquad
w(x,t):=\frac{v(x_{0}+Rx,t_{0}+K^{-\kappa}R^{\theta}t)}{KR^{\alpha}}
\quad\text{in }Q_{1}^{\theta}.
\end{align*}
Then $0<w\leq1$ in $Q_{1}^{\theta}$ and $w$ satisfies
\begin{align*}
|Dw|^{\kappa}\widetilde{F}(D^{2}w)-w_{t}
&=K^{-(\kappa+2-\gamma)}g_{\epsilon/(KR^{\alpha})}(w) \quad \text{in } Q_{1}^{\theta},
\end{align*}
in the viscosity sense, where $\widetilde{F}$ is defined by
\begin{align*}
\widetilde{F}(X)
:=\frac{F(KR^{\alpha-2}X)}{KR^{\alpha-2}} \quad\text{for }X\in\mathcal{S}^{n},
\end{align*}
and satisfies \ref{a1} and \ref{a2} with the same ellipticity constants.
Since $0<K^{-(\kappa+2-\gamma)}\leq\delta_{0}$, applying Lemma \ref{iterated growth} and scaling back, we obtain
\begin{align}
\label{adapted growth}
\sup_{B_{r}(x_{0})\times(t_{0}-K^{-\kappa}r^{\theta},t_{0}]}v
\leq C\bigl(v(z_{0})+Kr^{\alpha}\bigr)
\qquad \text{for } 0<r\leq R.
\end{align}
If $K^{\kappa/\theta}r\leq R$, apply \eqref{adapted growth} with
$K^{\kappa/\theta}r$ in place of $r$ to obtain
\begin{align*}
\sup_{Q_{r}^{\theta}(z_{0})}v
\leq \sup_{B_{K^{\kappa/\theta}r}(x_{0})\times(t_{0}-r^{\theta},t_{0}]} v
\leq C\left(v(z_{0})
+K^{1+\frac{\alpha\kappa}{\theta}}r^{\alpha}\right)=C(v(z_{0})+K^{\frac{\kappa+2-\gamma}{2-\gamma}}r^{\alpha}),
\end{align*}
where we used
\begin{align*}
1+\frac{\alpha\kappa}{\theta}
=\frac{\kappa+2-\gamma}{2-\gamma}.
\end{align*}
If $K^{\kappa/\theta}r>R$, then $K\geq M/R^{\alpha}$ gives
\begin{align*}
\sup_{Q_{r}^{\theta}(z_{0})}v\leq M\leq KR^{\alpha}<K^{1+\frac{\alpha\kappa}{\theta}}r^{\alpha}
=K^{\frac{\kappa+2-\gamma}{2-\gamma}}r^{\alpha}.
\end{align*}
Finally, the definition of $K$ gives
\begin{align}
\label{k upper}
K^{\frac{\kappa+2-\gamma}{2-\gamma}}=\max\left\{1,\delta_{0}^{-1/(2-\gamma)},
\left(\frac{M}{R^{\alpha}}\right)^{
\frac{\kappa+2-\gamma}{2-\gamma}}\right\}
\leq C\left[
1+\left(\frac{M}{R^{\alpha}}\right)^{
\frac{\kappa+2-\gamma}{2-\gamma}}
\right].
\end{align}
Combining these estimates proves \eqref{optimal growth}.
\end{proof}

\subsection{Forward growth}
To extend the spatial bound at $t=t_0$ to later times, we consider the profile $(|x-x_0|^\theta+c(t-t_0))^{1/(2-\gamma)}$, which reduces to $|x-x_0|^\alpha$ at $t=t_0$.
We then choose $c>0$ so that this function is a supersolution.

\begin{prop}
\label{forward growth proposition}
Let $F$ satisfy \ref{a1} and \ref{a2}. 
Let $v$ be a positive viscosity solution to \eqref{pe2 eq} in $B_{R}(x_{0})\times(t_{0}-R^{\theta},t_{0}+R^{\theta})$, 
with $\epsilon>0$ and $\delta=1$. Assume that $0<v\leq M$ in
$B_{R}(x_{0})\times(t_{0}-R^{\theta},t_{0}+R^{\theta})$. 
Then
\begin{align}
\label{forward growth}
\sup_{B_{r}(x_{0})\times[t_{0},t_{0}+r^{\theta})}v
\leq C\left\{v(x_{0},t_{0})+
\left[1+\left(\frac{M}{R^{\alpha}}\right)^{\frac{\kappa+2-\gamma}{2-\gamma}}\right]r^{\alpha}\right\} \qquad \text{for } 0<r\leq R,
\end{align}
where $C=C(n,\lambda,\Lambda,\kappa,\gamma)$.
\end{prop}

\begin{proof}
Write $z_{0}=(x_{0},t_{0})$. After translation, assume $z_{0}=(0,0)$,
and choose $K$ as in the proof of Theorem~\ref{optimal growth proposition}.
By \eqref{adapted growth}, at $t=0$ we have
\begin{align}
\label{initial}
v(x,0)\leq C_{2}v(0,0)+C_{2}K|x|^{\alpha}
\quad\text{for }x\in B_{R},
\end{align}
where $C_{2}=C_{2}(n, \lambda, \Lambda, \kappa, \gamma)\geq1$.
Set
\begin{align*}
A:=C_{2}K,\quad
B:=\Lambda\left(\frac{\theta}{2-\gamma}\right)^{\kappa}
\theta(n+\theta-2), \quad
W(x,t):=C_{2}v(0,0)
+A\bigl(|x|^{\theta}+BA^{\kappa}t\bigr)^{1/(2-\gamma)}.
\end{align*}
For $x\in B_{R}\setminus\{0\}$ and $t\in(0,R^{\theta})$,
write $S=|x|^{\theta}+BA^{\kappa}t$.
Direct calculation gives
\begin{align}
\label{DC}
|DW(x,t)|=\frac{A\theta}{2-\gamma}
S^{\frac{\gamma-1}{2-\gamma}}|x|^{\theta-1}, \quad
W_{t}(x,t)=\frac{BA^{\kappa+1}}{2-\gamma}S^{\frac{\gamma-1}{2-\gamma}}.
\end{align}
In addition, by Lemma~\ref{Hessian}, the eigenvalues of $D^{2}W(x,t)$ are
\begin{align*}
\frac{A\theta}{2-\gamma}S^{\frac{\gamma-1}{2-\gamma}}|x|^{\theta-2}>0
\end{align*}
in the $n-1$ tangential directions and
\begin{align*}
\frac{A\theta}{2-\gamma}S^{\frac{\gamma-1}{2-\gamma}}
|x|^{\theta-2}
\left[\theta-1+\frac{\theta(\gamma-1)}{2-\gamma}
\frac{|x|^{\theta}}{S}\right]>0
\end{align*}
in the radial direction, where we used
\begin{align*}
\theta-1+\frac{\theta(\gamma-1)}{2-\gamma}\frac{|x|^{\theta}}{S}
\geq\theta-1+\frac{\theta(\gamma-1)}{2-\gamma}
=\frac{\theta}{2-\gamma}-1=\alpha-1>0.
\end{align*}
Also note that
\begin{align}
\label{identity}
\kappa(\theta-1)+\theta-2
=\frac{\kappa(\kappa+2)(1-\gamma)}{\kappa+2-\gamma}>0, \quad
\frac{(\gamma-1)(\kappa+1)}{2-\gamma}
+\frac{\kappa(\theta-1)+\theta-2}{\theta}
=\frac{\gamma-1}{2-\gamma}.
\end{align}
Thus, by \ref{a1}, \eqref{DC}, $|x|^{\theta}\leq S$, and \eqref{identity}, we obtain
\begin{align*}
|DW|^{\kappa}F(D^{2}W)
&\leq |DW|^{\kappa}\M^{+}(D^{2}W) \\
&\leq\Lambda\left(\frac{A\theta}{2-\gamma}\right)^{\kappa+1}
(n+\theta-2)S^{\frac{(\gamma-1)(\kappa+1)}{2-\gamma}}
|x|^{\kappa(\theta-1)+\theta-2}\\
&\leq\Lambda\left(\frac{A\theta}{2-\gamma}\right)^{\kappa+1}(n+\theta-2)
S^{\frac{(\gamma-1)(\kappa+1)}{2-\gamma}
+\frac{\kappa(\theta-1)+\theta-2}{\theta}}\\
&=\Lambda\left(\frac{A\theta}{2-\gamma}\right)^{\kappa+1}(n+\theta-2)
S^{\frac{\gamma-1}{2-\gamma}}=W_{t}
\end{align*}
in $(B_{R}\setminus\{0\})\times(0,R^{\theta})$.
We now assume that $\xi\in C^{2,1}$ touches $W$ from below at
$(0,t)$ with $0<t<R^{\theta}$. Since $W$ is differentiable there,
we have $D\xi(0,t)=DW(0,t)=0$ and
$\xi_{t}(0,t)=W_{t}(0,t)>0$. Hence,
\begin{align*}
|D\xi(0,t)|^{\kappa}F(D^{2}\xi(0,t))-\xi_{t}(0,t)
=-W_{t}(0,t)<0\leq g_{\epsilon}(W(0,t)).
\end{align*}
Therefore, $W$ is a viscosity supersolution
to \eqref{pe2 eq} with $\delta=1$ in $B_{R}\times(0,R^{\theta})$.
We now extend $v$ to $\overline{B_R}\times[0,R^\theta]$ by taking its upper semicontinuous envelope, retaining the same notation. 
The extension is bounded above by $M$ and agrees with $v$ in $B_R\times[0,R^\theta)$.
Using \eqref{initial} and the definitions of $A$ and $K$, we have
\begin{align*}
v\leq W \quad\text{on }\overline{B_{R}}\times\{0\}, \quad
v\leq M\leq KR^{\alpha}\leq AR^{\alpha}\leq W
\quad\text{on }\partial B_{R}\times[0,R^{\theta}].
\end{align*}
Since $v$ is a viscosity subsolution to \eqref{pe2 eq} with $\delta=1$,
Lemma~\ref{comparison} gives $v\leq W$ in
$B_{R}\times[0,R^{\theta})$. For $|x|< r$ and
$0\leq t<r^{\theta}$, we conclude that
\begin{align*}
v(x,t)
\leq C_{2}v(0,0)+A(1+BA^{\kappa})^{1/(2-\gamma)}r^{\alpha}
\leq C\left(v(0,0)+K^{\frac{\kappa+2-\gamma}{2-\gamma}}r^{\alpha}\right).
\end{align*}
Finally, \eqref{k upper} proves \eqref{forward growth}.
\end{proof}

We now complete the proof of Theorem~\ref{thm1} by passing to the limit in the preceding growth estimates.
\begin{proof}[Proof of Theorem \ref{thm1}]
Proposition~\ref{existence and equiconti} gives a nonnegative
viscosity solution $u\in C(\overline{\Omega_{T}})$ to \eqref{eq} and a
sequence $\epsilon_{j}\to0^{+}$ such that
$u_{\epsilon_{j}}\to u$ locally uniformly. Fix $z_{0}=(x_{0},t_{0})$ with
$u(z_{0})=0$. Apply \eqref{optimal growth} with
$M=\|\varphi\|_{L^{\infty}(\partial_{p}\Omega_{T})}+\epsilon_{j}$ and let $j\to\infty$.
Since $u_{\epsilon_{j}}(z_{0})\to0$, we obtain
\begin{align*}
\sup_{Q_{r}^{\theta}(z_{0})}u
\leq C\left[1+
\left(\frac{\|\varphi\|_{L^{\infty}(\partial_{p}\Omega_{T})}}{R^{\alpha}}\right)^{\frac{\kappa+2-\gamma}{2-\gamma}}\right]r^{\alpha}.
\end{align*}
Passing to the limit in \eqref{forward growth} proves the
corresponding estimate in $B_{R}(x_{0})\times
 (t_{0}-R^{\theta},t_{0}+R^{\theta})\Subset\Omega_{T}$.
\end{proof}

\section{Non-degeneracy}
\label{sec6}   
We prove Theorem~\ref{thm2} directly from the viscosity
subsolution inequality under assumption~\ref{a1} alone,
without using the approximation procedure.

\begin{proof}[Proof of Theorem \ref{thm2}]
By continuity, it is enough to prove the result for $z_{0}\in\left\{u>0\right\}\cap \Omega_{T}$ with $Q_{r}^{\theta}(z_{0})\Subset \Omega_{T}$.
After translation,
assume that $z_{0}=(0,0)$, and define
\begin{align*}
M:=\max_{\overline{Q_{r}^{\theta}}}u>0,\qquad
\Psi_{\eta}(x,t):=M\left(\frac{|x|^{2}}{r^{2}}
-\frac{t}{r^{\theta}}\right)+\frac{\eta}{-t}
\quad(-r^{\theta}\leq t<0).
\end{align*}
By continuity, choose $s\in(-r^{\theta},0)$ sufficiently close to zero and then $\eta>0$ sufficiently small so that
\begin{align}
\label{positive pt}
u(0,s)-\Psi_{\eta}(0,s)=u(0,s)+\frac{Ms}{r^{\theta}}+\frac{\eta}{s}>0.
\end{align}
Choose $h\in(0,-s)$ such that
$\eta/h>M$. 
Since $\Psi_{\eta}(x,t)\geq M\geq u$ when $|x|=r$ or $t\in\left\{-r^{\theta},-h\right\}$, we have
\begin{align*}
u-\Psi_{\eta}\leq0 \quad \text{on } (\ov{B_{r}}\times\{-r^{\theta},-h\})\cup (\partial B_{r}\times[-r^{\theta},-h]).   
\end{align*}
Thus, by \eqref{positive pt}, $u-\Psi_{\eta}$ attains a positive maximum at some
$(x_{*},t_{*})\in B_{r}\times(-r^{\theta},-h)$, and
$u(x_{*},t_{*})>\Psi_{\eta}(x_{*},t_{*})>0$.
Direct differentiation gives
\begin{align*}
D\Psi_{\eta}=\frac{2M}{r^{2}}x,\qquad
D^{2}\Psi_{\eta}=\frac{2M}{r^{2}}I,\qquad
(\Psi_{\eta})_{t}=-\frac{M}{r^{\theta}}
+\frac{\eta}{(-t)^{2}}.
\end{align*}
Since $\gamma-1<0$, the viscosity subsolution inequality and
\ref{a1} yield
\begin{align*}
M^{\gamma-1}
&\leq u(x_{*},t_{*})^{\gamma-1}\\
&\leq |D\Psi_{\eta}(x_{*},t_{*})|^{\kappa}F(D^{2}\Psi_{\eta}(x_{*},t_{*}))
-(\Psi_{\eta})_{t}(x_{*},t_{*})\\
&\leq2^{\kappa+1}n\Lambda\frac{M^{\kappa+1}}{r^{\kappa+2}}
+\frac{M}{r^{\theta}}-\frac{\eta}{(-t_{*})^{2}}\\
&\leq2^{\kappa+1}n\Lambda\frac{M^{\kappa+1}}{r^{\kappa+2}}
+\frac{M}{r^{\theta}}.
\end{align*}
Dividing by $M^{\gamma-1}$ and using
$\alpha(\kappa+2-\gamma)=\kappa+2$ and $\theta=\alpha(2-\gamma)$,
we obtain
\begin{align*}
1\leq2^{\kappa+1}n\Lambda
\left(\frac{M}{r^{\alpha}}\right)^{\kappa+2-\gamma}
+\left(\frac{M}{r^{\alpha}}\right)^{2-\gamma}.
\end{align*}
Choose $c_{0}=c_{0}(n, \Lambda, \kappa, \gamma)>0$ so that
\begin{align*}
2^{\kappa+1}n\Lambda c_{0}^{\kappa+2-\gamma}
+c_{0}^{2-\gamma}\leq1.
\end{align*}
The preceding inequality implies $M\geq c_{0}r^{\alpha}$.
\end{proof}

\begin{rmk}
An analogous lower bound on a time slice or in a forward cylinder
may fail at a free boundary point. Indeed,
\begin{align*}
u(x,t)=\bigl[(2-\gamma)(t_{*}-t)\bigr]_{+}^{1/(2-\gamma)}
\end{align*}
solves the equation in its positive set and vanishes identically for $t\geq t_{*}$.  
In particular,
\begin{align*}
\sup_{B_{r}(x_{0})}u(\cdot,t_{*})
=\sup_{B_{r}(x_{0})\times[t_{*},t_{*}+r^{\theta}]}u=0.
\end{align*}
\end{rmk}

\vspace{0.2cm}
{\bf Acknowledgments.} This publication is based upon work supported by King Abdullah University of Science and Technology (KAUST) under Award No. ORFS-CRG12-2024-6430.

\vspace{0.2cm}
{\bf AI Disclosure.} AI tools were used to improve the writing, assist with reference searches, and help check the manuscript for errors and inconsistencies.

\vspace{0.2cm}
{\bf Data Availability.}
Data sharing does not apply to this article as no datasets were generated or analyzed during the current study.

\providecommand{\bysame}{\leavevmode\hbox to3em{\hrulefill}\thinspace}
\providecommand{\MR}{\relax\ifhmode\unskip\space\fi MR }
\providecommand{\MRhref}[2]{%
  \href{http://www.ams.org/mathscinet-getitem?mr=#1}{#2}
}
\providecommand{\href}[2]{#2}


\begin{thebibliography}{10}

\bibitem{AC1981}
Hans~W. Alt and Luis~A. Caffarelli, \emph{Existence and regularity for a
  minimum problem with free boundary}, J. Reine Angew. Math. \textbf{325}
  (1981), 105--144. \MR{618549}

\bibitem{AP1986}
Hans~W. Alt and Daniel Phillips, \emph{A free boundary problem for semilinear
  elliptic equations}, J. Reine Angew. Math. \textbf{368} (1986), 63--107.
  \MR{850615}

\bibitem{AS2022}
Damião~J. Araújo and Ginaldo~S. S\'a, \emph{Infinity {L}aplacian equations
  with singular absorptions}, Calc. Var. Partial Differential Equations
  \textbf{61} (2022), no.~4, Paper No. 132, 16. \MR{4417397}

\bibitem{ASU2024}
Damião~J. Araújo, Ginaldo~S. S\'a, and José~Miguel Urbano, \emph{Sharp
  regularity for a singular fully nonlinear parabolic free boundary problem},
  J. Differential Equations \textbf{389} (2024), 90--113. \MR{4693222}

\bibitem{AT2013}
Damião~J. Araújo and Eduardo~V. Teixeira, \emph{Geometric approach to
  nonvariational singular elliptic equations}, Arch. Ration. Mech. Anal.
  \textbf{209} (2013), no.~3, 1019--1054. \MR{3067831}

\bibitem{ASP26}
Alessandro Audrito and Tom\'as Sanz-Perela, \emph{On the existence of
  solutions to some singular parabolic free boundary problems}, Nonlinear
  Anal. \textbf{272} (2026), Paper No. 114170.

\bibitem{BK26}
Sun-Sig Byun and SeungHyun Kim, \emph{On fully nonlinear free boundary problem
  with logarithmic singular absorption terms}, Calc. Var. Partial Differential
  Equations \textbf{65} (2026), no.~3, Paper No. 82. \MR{5033048}

\bibitem{C1977}
Luis~A. Caffarelli, \emph{The regularity of free boundaries in higher
  dimensions}, Acta Math. \textbf{139} (1977), no.~3-4, 155--184. \MR{454350}

\bibitem{CIL1992}
Michael~G. Crandall, Hitoshi Ishii, and Pierre-Louis Lions, \emph{User's guide
  to viscosity solutions of second order partial differential equations}, Bull.
  Amer. Math. Soc. (N.S.) \textbf{27} (1992), no.~1, 1--67. \MR{1118699}

\bibitem{DD16}
Nguyen~Anh Dao and Jesus~Ildefonso D\'iaz, \emph{A gradient estimate to a
  degenerate parabolic equation with a singular absorption term: the global
  quenching phenomena}, J. Math. Anal. Appl. \textbf{437} (2016), no.~1,
  445--473. \MR{3451975}

\bibitem{DM05}
Juan D\'avila and Marcelo Montenegro, \emph{Existence and asymptotic behavior
  for a singular parabolic equation}, Trans. Amer. Math. Soc.
  \textbf{357} (2005), no.~5, 1801--1828. \MR{2115077}

\bibitem{DS21}
Daniela De~Silva and Ovidiu Savin, \emph{On certain degenerate one-phase free
  boundary problems}, SIAM J. Math. Anal. \textbf{53} (2021), no.~1, 649--680.
  \MR{4205089}

\bibitem{KU26}
Seunghyun Kim and Jos{\'e}~Miguel Urbano, \emph{On the fully nonlinear
  quenching problem}, Ann. Sc. Norm. Super. Pisa Cl. Sci. (5), to appear,
  \href{https://doi.org/10.2422/2036-2145.202509_013}{doi:10.2422/2036-2145.202509\_013}.

\bibitem{LLY2024}
Ki-Ahm Lee, Se-Chan Lee, and Hyungsung Yun, \emph{{$C^{1,\alpha}$}-regularity
  for solutions of degenerate/singular fully nonlinear parabolic equations}, J.
  Math. Pures Appl. (9) \textbf{181} (2024), 152--189. \MR{4676644}

\bibitem{LS2001}
Ki-Ahm Lee and Henrik Shahgholian, \emph{Regularity of a free boundary for
  viscosity solutions of nonlinear elliptic equations}, Comm. Pure Appl. Math.
  \textbf{54} (2001), no.~1, 43--56. \MR{1787106}

\bibitem{LLYZ25}
Se-Chan Lee, Yuanyuan Lian, Hyungsung Yun, and Kai Zhang, \emph{Time derivative
  estimates for parabolic $p$-{Laplace} equations and applications to optimal
  regularity}, \href{https://arxiv.org/abs/2503.04384}{arXiv:2503.04384}, 2025.

\bibitem{D1983-2}
Daniel Phillips, \emph{Hausdorff measure estimates of a free boundary for a
  minimum problem}, Comm. Partial Differential Equations \textbf{8} (1983),
  no.~13, 1409--1454. \MR{714047}

\bibitem{D1983}
\bysame, \emph{A minimization problem and the regularity of solutions in the
  presence of a free boundary}, Indiana Univ. Math. J. \textbf{32} (1983),
  no.~1, 1--17. \MR{684751}

\bibitem{RT2011}
Gleydson~C. Ricarte and Eduardo~V. Teixeira, \emph{Fully nonlinear singularly
  perturbed equations and asymptotic free boundaries}, J. Funct. Anal.
  \textbf{261} (2011), no.~6, 1624--1673. \MR{2813483}

\bibitem{T2018}
Eduardo~V. Teixeira, \emph{Nonlinear elliptic equations with mixed
  singularities}, Potential Anal. \textbf{48} (2018), no.~3, 325--335.
  \MR{3779091}

\bibitem{W1992}
Lihe Wang, \emph{On the regularity theory of fully nonlinear parabolic
  equations. II}, Comm. Pure Appl. Math. \textbf{45} (1992), no.~2,
  141--178. \MR{1139064}  

\bibitem{W1998}
Georg~S. Weiss, \emph{Partial regularity for weak solutions of an elliptic free
  boundary problem}, Comm. Partial Differential Equations \textbf{23} (1998),
  no.~3-4, 439--455. \MR{1620644}

\bibitem{W07}
Michael Winkler, \emph{Nonuniqueness in the quenching problem},
Math. Ann. \textbf{339} (2007), 559--597.  

\bibitem{WY22}
Yijing Wu and Hui Yu, \emph{On the fully nonlinear {A}lt-{P}hillips equation},
  Int. Math. Res. Not. IMRN (2022), no.~11, 8540--8570. \MR{4425844}

\end{thebibliography}
\end{document}